\documentclass{article}
\usepackage{amssymb,amsmath,dsfont}
\usepackage[latin1]{inputenc}
\usepackage[usenames]{color}

\def\Om       {\Omega}
\def\eps	{\varepsilon}

\def\v        {{\boldsymbol v}}

\def\u        {{\boldsymbol u}}
\def\w        {{\boldsymbol w}}

\def\x        {{\boldsymbol x}}

\def\b        {{\boldsymbol b}}

\def\f        {{\boldsymbol f}}
\def\C	{{\boldsymbol C}}
\def\H        {{\boldsymbol H}}

\def\L        {{\boldsymbol L}}
\def\X        {{\boldsymbol X}}
\def\G        {{\boldsymbol G}}
\def\J        {{\boldsymbol J}}
\def\n        {{\boldsymbol n}}

\def\fra#1#2{\displaystyle\frac{\mathstrut #1}{\mathstrut #2}}
\def\jnt      {\displaystyle\int}

\newtheorem{theorem}{Theorem}
\newtheorem{lemma}{Lemma}
\newtheorem{proposition}{Proposition}

\newtheorem{definition}{Definition}

\newtheorem{remark}{Remark}
\begin{document}
\title{An Analysis of the Evolution Equations for a Generalized Bioconvective Flow}

\author{ J.L. Boldrini\thanks{FEM-UNICAMP, Campinas, SP, Brazil. E-mail:
{\tt josephbold@gmail.com}. }, 
 B. Climent-Ezquerra\thanks{ Dpto. EDAN, Universidad de Sevilla, Spain. E-mail:
{\tt bcliment@us.es}. Partially supported by Ministerio de Ciencia, Innovaci\'on y Universidades, Grant PGC2018-098308-B-I00, with the collab oration of FEDER. },
 M.A. Rojas-Medar\thanks{Departamento de Matem\'atica, Universidad de Tarapac\'a, Arica, Chile. E-mail: {\tt marko.medar@gmail.com}. Partially supported by MATH-AMSUD project 21-MATH-03 (CTMicrAAPDEs), and CAPES-PRINT 88887.311962/2018-00 (Brazil).},   
  M.D. Rojas-Medar\thanks{Departamento de
Matem\'aticas, Facultad de Ciencias B\'asica, Universidad de Antofagasta,
Antofagasta, Chile. E-mail: {\tt mrojas@uantof.cl}. Partially supported by COL\'OQUIOS DE MATEM\'ATICA of Universidad de Antofagasta.}
}

\date{}
\maketitle
\begin{abstract}
We prove results on existence and uniqueness of solutions of a system of equations modeling the evolution of a generalized bioconvective flow. 
The mathematical model considered in the present work describes the convective motion generated by the upward swimming of a culture of microorganisms under the influence of vertical gravitational forces, in an incompressible viscous fluid whose viscosity may  depend on the concentration of microorganisms. 
\end{abstract}
Keywords: Bioconvective flow, weak solution, strong solution.

\section{Introduction}
\label{sec:1}
The purpose of this paper is to study the existence and uniqueness of evolution solutions of a system of equations modeling a generalized bioconvective flow.

Bioconvective flows are important biological phenomena associated to the appearance of motion in bodies of fluid where there are cultures of microorganisms that tend to swim upward against a vertical gravitational force.
This swimming leads to differences of microorganism concentration, causing convective flows in the fluid, which in turn change the microorganism concentration.

Simple mathematical modeling of such situation assumes constant fluid viscosity; 
however, this simplification is not adequate in general since the presence of microorganisms implies in an effective flow  viscosity that may depend on the microorganism concentration.
Thus, more realistic models must consider the effects of such concentration on the fluid viscosity;
this leads to the generalized bioconvective models, whose associated fluid motions are called generalized bioconvective flows. 

In the present work, we consider the following generalized bioconvective model:
\begin{equation} 
\label{sistema1}
\left\{\begin{array}{l} \displaystyle\frac{\partial \u}{\partial t}
- 2 \ {\rm div}\, (\mu(c)D(\u))
+ \u\cdot  \nabla \u + \nabla p = -g(1+\rho c)\chi+\f, \\
{\rm div}\, \u = 0, \\
\displaystyle\frac{\partial c}{\partial t} - \theta \Delta c +
\u\cdot \nabla c + U \displaystyle\frac{\partial c}{\partial x_3} =
0, \qquad {\rm in} \quad (0, T) \times \Omega.
\end{array}\right.
\end{equation}

\noindent
Here, the following notation is used:
$\Omega \subset \mathbb{R}^3$ is a bounded domain representing the region where the fluid flow occurs.
The unknowns, $\u(x,t)$, $p(x,t)$ and $c(x,t)$ represent the fluid velocity, the hydrostatic pressure and the concentration of microorganisms, respectively, at a point
$x = (x_1, x_2, x_3) \in \Omega$, at time $t \in~[0,T]$, here $0<T\leq +\infty$ is the final time of interest.
Function $\mu(\cdot) > 0$ gives the kinematic viscosity of fluid.
The positive constant $g$ is the  intensity of the acceleration of gravity,
and function $\f (\x,t)$ is a given external force.
Parameters $\theta$,  $\chi$, $U$ and $\rho$ are assumed to be constants. 
Parameter $\theta$ is positive and gives the rate of diffusion of microorganisms;
$\chi$ is a unitary vector in the vertical direction, i.e. $\chi = (0, 0, 1)^t$, and it is placed to represent that the gravitational forces acting
along the vertical; 
$U$ is  positive  and denotes the average velocity of upward swimming of  microorganisms;
$\rho$ is positive and is defined by $\rho = \displaystyle\frac{\rho_0}{\rho_m} - 1$, where $\rho_0$
and $\rho_m$ are the density of one organism and the culture fluid
density, respectively.

The symbols $\nabla, \Delta$ and ${\rm div}$ represent
the gradient, Laplacian and divergence operator, respectively;
operator
$
D(\u) =\displaystyle\frac12 (\nabla \u + (\nabla \u)^t)
$
is the symmetrized gradient of velocity.

In system $(\ref{sistema1})$, the first two equations model an incompressible Navier-Stokes with concentration dependent viscosity and a right-hand side where the influence of the concentration of microorganism is taken in consideration;
the third  equation is a conservation equation of form
$\displaystyle\frac{\partial c}{\partial t} + {\rm div} \ \mathcal{F} = 0, \ x \in \Omega, \ t > 0$,
with $\mathcal{F} = \u c  - \theta \nabla c + Uc \chi$. 
For the derivation of system (\ref{sistema1}), see Levandowsky {\it et al.}  \cite{levan}.

In order to simplify the notation, as in previous works, we set: 
\[
c = (g\rho)^{-1} m, \quad p=q-gx_3, \quad \mu(c) = \nu(m) ,
\]

\noindent
and rewrite system (\ref{sistema1}) as
\begin{equation} 
\label{sistema2}
\left\{\begin{array}{l} \displaystyle\frac{\partial \u}{\partial t}
- 2 \ {\rm div}\, (\nu(m)D(\u))
+ \u\cdot \nabla \u + \nabla q = -m \chi +\f, \\
{\rm div}\, \u = 0, \\
\displaystyle\frac{\partial m}{\partial t} - \theta \Delta m +
\u\cdot \nabla m + U \displaystyle\frac{\partial m}{\partial x_3} =
0, \qquad {\rm in} \quad (0, T) \times \Omega.
\end{array}\right.
\end{equation}

We will analyze system (\ref{sistema2})  together with initial conditions:
\begin{equation}
\label{ic}
\u = \u_0 \quad \mbox{and} \quad m = m_0 \quad \mbox{at time} \, t =0 ,
\end{equation}

\noindent
and the following boundary conditions:
\begin{equation}
\label{bc0}
\left\{\begin{array}{l}
\u = 0 \qquad {\rm on} \qquad (0,T] \times S, \\
\u\cdot \n = 0 \qquad {\rm on} \qquad (0, T] \times \Gamma,\\
\nu(m)[D(\u)\n - \n\cdot(D(\u)\n)\n] =\b_1
\qquad {\rm on} \qquad (0, T] \times \Gamma, 
\vspace{0.1cm}
\\
\theta \displaystyle\frac{\partial m}{\partial \n} - Umn_3 = 0
\qquad {\rm on} \qquad  (0, T] \times \partial\Omega. 
\end{array}\right.
\end{equation}

The boundary $\partial\Omega = S \cup \overline{\Gamma} $  consists of a rigid nonslip part, $S$, and a free surface, $\Gamma$;
we assume that $S \neq \emptyset$ and  that $\overline{S}$ and disjoint $\overline{\Gamma}$.

The vector 
$\n(x)=(n_{1}(x), n_{2}(x), n_{3}(x))$ is the exterior unitary normal vector at the point $x \in
\partial \Omega$ and $\displaystyle\frac{\partial}{\partial
\n}$ is the normal derivative on $\partial \Omega$. The second and third equations of (\ref{bc0}) represent the usual slip conditions for the fluid velocity
(see Joseph \cite{joseph}). In particular, the third one means that the tangential component of the stress vector has a value $\b_1(x)$ verifying $\b_1(x)\cdot\n(x)=0$ on $(0, T] \times \Gamma$. 
The last boundary condition in (\ref{bc0}) 
establishes the condition of null microorganism flow at each point $x \in
\partial\Omega$. 
In fact, the given boundary conditions for the velocity imply that 
the normal component of the microorganism flux at the boundary, written in terms of $c$, becomes
$\mathcal{F} \cdot \n = (\u c  - \theta \nabla c + Uc \chi) \cdot \n= \theta \displaystyle\frac{\partial c}{\partial \n} - U c n_3 = 0$
on $(0, T] \times \partial\Omega$; this, rewritten in terms of $m$ becomes exactly the last equation in (\ref{bc0}).

\vspace{0.1cm}
We finally observe that the last equation equation in \eqref{sistema2}  and the boundary conditions imply that
\begin{equation}
\label{TotalMassConservation}
\int_\Omega m (x, t) dx = \int_\Omega  m_0 (x) dx = \alpha, \quad \mbox{for all} \; t \in [0,T] ,
\end{equation}

\vspace{0.2cm}

Now, let us briefly comment on previous published works directly related to the present one.

The particular case of bioconvective flow with constant viscosity was introduced independently in Levandowsky {\it et al.}  \cite{levan} and Moribe \cite{moribe};
by using  biological and physical arguments these authors present qualitative results. 
Rigorous mathematical analysis of this classical case of constant viscosity $\mu$ and Dirichlet boundary conditions for the velocity
 was carried out in Kan-On {\it et al.} \cite{kan}.
In \cite{kan}, the authors prove results about existence of solutions for both the stationary and evolution evolution problems;
they use the Galerkin method for proving the existence of weak solutions and the semigroup approach,  together with the method of successive approximations, for proving the existence of strong solutions.

Concerning the boundary conditions for the velocity, Solonnikov and Scadilov \cite{solonnikov} and also Mulone and Salemi \cite{mulone}
considered the same type of boundary conditions presented in this work. 
However, the authors of \cite{solonnikov} and  \cite{mulone} considered just the classical Navier-Stokes equations.

As for bioconvective models dealing with non-constant the viscosities,
Boldrini {\it et al.}  \cite{BRMRM1} studied the existence of stationary solutions, while  Climent-Ezquerra {\it et al.} \cite{CFR}
analyzed problems with time-periodic conditions.

\vspace{0.1cm}
As we have already said, we are interested here in evolution solutions of the generalized case;
that is, when the concentration of microorganisms affects the effective viscosity of the fluid.
A consequence is that the term ${\rm div}\,(\mu(c)D(\u))$ brings a much higher degree of difficulty for the analysis of the problem
as compared to the classical case with constant $\mu$. 
The adaptation of techniques and arguments used in Lorca and Boldrini \cite{lorca1}, \cite{lorca2} for a generalized Boussinesq model will be
fundamental to our technical arguments. 
We also remark that, in previous works related to this one, the proof of existence of strong solutions required in an  essential way the condition that $\b_1=0$; this is also the case of the present paper. 

\vspace{0.2cm}
In summary, the main novelties in this paper are the following:
\begin{itemize}
\item We consider a viscosity term depending on the concentration of microorganisms  
and also a slip boundary condition on part of the boundary. These features introduce significant difficulties in the analysis.

\item We obtain weak solution of the problem without imposing smallness conditions.

\item We prove the existence of local and global strong solution by using a spectral Galerkin method instead of using semigroups.
 The derived estimates for the spectral Galerkin approximations could be used to derive error estimates.
 
 \item Unlike \cite{kan}, we obtain weak and strong solutions without reduction to decay problem.

\end{itemize}

\vspace{0.1cm}
The organization of the paper is the following: 

In Section  \ref{NotationFunctionalSpacesMainResults}, we fix the notation, recall certain functional spaces and 
state several auxiliary technical results that will be used throughout the paper.
Section \ref{ExistenceOfGlobalEvolutionaryWeakSolutions} is devoted to the proof of the global existence of weak evolution solutions.
Section  \ref{sec:4} deals with the existence and uniqueness of local strong evolution solutions.
Section \ref{GlobalEvolutionaryStrongSolutions} is devoted to the proof of existence of global strong evolution solutions. 
Finally, in Section \ref{ExponentialL2StabilityOfSmallStrongStationarySolutions},
we prove the exponential $L^2$-stability of small strong 

\section{Notation, functional spaces and some technical lemmas}
\label{NotationFunctionalSpacesMainResults}

\subsection{Notation and some interpolation results}

Let $\Omega \subseteq \mathbb{R}^3$ be a bounded domain with boundary of class $C^2$. 
The $L^p$-norm on $\Omega$ is denoted by $\mid \cdot \mid_p$, $1\leq p \leq \infty$; 
the usual inner product in $L^2 (\Omega)$ is denoted by $(\cdot,\cdot)$,
and, to simplify the notation, its  associated $L^2$-norm, just by $| \cdot |$ instead of $| \cdot |_2$. 
Let $H^m(\Omega)$ be the usual Sobolev spaces on $\Omega$ with norm $\| \cdot \|_m$. 
Moreover, $H^1_0(\Omega)$ denotes the completion of $C^\infty _0(\Omega)$ under the norm $\parallel\cdot\parallel_1$.
The generalization of such Sobolev spaces are denoted $H^\alpha (\Omega)$, $\alpha \in \mathbb{R}$, the space of Bessel potentials,
with corresponding norm denoted by $\| \cdot \|_\alpha$.

If $B$ is a Banach space,  $L^q(0,T;B)$ is the Banach space of the
$B$-valued functions defined in the interval $(0,T)$ that are $L^q$-integrable
in the sense of Bochner;
the norms in such spaces are denoted $\| \cdot \|_{L^q (0,T; B)}$.

Boldface letters will be used for spaces associated to multiple copies of the same space; for instance, $\L^2 (\Omega) =L^2(\Omega)^N$.

\begin{remark}
As usual in works dealing with partial differential equations, 
in the derivation of the estimates,
we will use a positive generic constant $C$ depending only of $\Omega$ and given fixed data of the problem.
Sometimes,  when we fill that it is convenient to distinguish certain constants, other letters or subscripts may be used as well.
\end{remark}

\begin{lemma}(Some interpolation inequalities Adams \cite{adams}) 
\label{ineq}
\[
| v |_6\leq C\| v\|_1, 
\quad
| v |_3 \leq | v |^{1/2} \Vert v\Vert_1^{1/2}, 
\quad
| v |_4 \leq | v |^{1/4}\| v \|_1^{3/4} \quad\forall \, v \in H^1 (\Omega),
\]
\[
| v |_\infty \leq C \| v \|_1^{1/2} \| v \|_2^{1/2}\quad\forall \, v \in H^2 (\Omega).
\]
\end{lemma}

\subsection{Specific functional spaces and some known technical results}

As we mentioned before, we assume that $\partial\Omega = S \cup \overline{\Gamma} $, with $\emptyset \neq S$ and  disjoint $\overline{S}$ and $\overline{\Gamma}$.

We define
\[
\dot \H  = \{\u \in \C^\infty_0 (\overline\Omega);\,
\u|_S = 0, \, \u \cdot \textbf{n}|_{\Gamma} = 0\},
\]

\noindent
with norm given by
\begin{equation} \label{norma}
||\u||_{\H}
= \left[\int_\Omega \nabla \u \cdot \nabla \u\, dx \,\right]^{1/2}
= [(\nabla \u, \nabla \u)]^{1/2} = |\nabla \u|.
\end{equation}
and 
\[
\H \; \mbox{is the clousure of}\; \dot \H(\Omega) \; \mbox{ with respect to the norm} \;  \|~ \|_{\H}.
\]

We denote by
$
\C^\infty_{0,\sigma} (\Omega)=
\{ \f \in \C^\infty_0 (\Omega) : \ {\rm div} \ \f = 0 \}
$
and
\[
\X \; \mbox{is the closure of $\C^\infty_{0,\sigma} (\Omega) $ in $\L^2 (\Omega)$}.
\]

 If $
\G =
\{\w \in \L^2 (\Omega), \w = \nabla q, q \in H^1 (\Omega)\}
$, it is well known that
$
\L^2 (\Omega) = \X \oplus \G
$ (see, for instance, Temam \cite{temam}). 

We define
$
\dot {\J}= \{\u \in \dot \H(\Omega), \
{\rm div} \ \u = 0\}
$

\noindent
and 
\[
\J_0 \; \mbox{is  the closure of $\dot \J$ in the norm $\|\cdot\|_{\H} $ defined in (\ref{norma})}.
\]

\vspace{0.2cm}
The computations to be done with the velocity will require the following technical results.

The proofs of the formulas presented in the next lemma can be seen in Boldrini, Rojas-Medar \cite{BRMRM1};
they are analogous to the ones given in Solonnikov, Scadilov \cite{solonnikov}.
\begin{lemma}
 \label{lema7} 
Let $\u, \v \in \J_{0} \,$
and $ q \in C^{1}$ then,
\[
\int_\Omega [-2{\rm div}\, (\nu(m)D(\u)) + \nabla q]\varphi \,dx =
2\int_\Omega\nu(m)D(\u):D(\varphi)\,dx- 2\int_{\gamma}b_{1}\cdot\varphi \, dS.
\]
Let $m, \phi \in B \,$ then,
\[
\displaystyle\int_\Omega- \theta \Delta m\, \phi\,dx + \int_\Omega U \displaystyle\frac{\partial m}{\partial x_3}\, \phi\,dx  
= 
\theta\int_\Omega  \nabla m\, \nabla\phi\,dx -U \int_\Omega  \displaystyle m\frac{\partial \phi}{\partial x_3}\,dx\]
\end{lemma}

\vspace{0.1cm}
The next two lemmas give important estimates.
\begin{lemma} 
\label{Lema 1.1} 
(Poincar\'e inequality). 
Let $\Omega$ be a bounded domain of $\mathbb{R}^3$ with boundary $\partial\Omega$, of class $C^1$,
and let $\emptyset \neq \Sigma \subseteq \partial\Omega$ be a part of the boundary.
Then, there exists a positive constant
$C_P$, depending only on  $\Omega$ and  $\Sigma$, such that
$$
|\u| \leq C_P |\nabla \u|, \quad \mbox{for all}
\quad \u \in (H^1(\Omega) )^3, \quad \mbox{such that} \quad \u |_\Sigma = 0.
$$
\end{lemma}
\begin{lemma} \label{Lema 1.2} (Korn inequality,  Solonnikov, Scadilov \cite{solonnikov}). Let $\Omega$ be a
bounded domain of $\mathbb{R}^3$
with boundary  $\partial \Omega$ of class $C^2$. Then, there exists a
positive constant $C$, such that:
$$
\|\u\|_{\H} = |\nabla \u| \leq C
|D(\u)|, \quad \forall \, \u \in \H.
$$
\end{lemma}

\noindent
Therefore, there exists a positive constant $C$  such that
\[
|\u|^2 \leq C |D(\u)|^2, \ \ \forall \, \u \in \H.
\]

\noindent
and  the norms $|\nabla \u|$ and $|D(\u)|$ are equivalent in $\H$.

\vspace{0.2cm}
Next, similarly as in Kan-on, Narukaw, Teramoto \cite{kan}  and Rionero, Mulone \cite{rionero}, we introduce two important linear operators for the analysis of strong solutions of the the problem considered in the present work.
\begin{lemma} 
\label{Lema 1.8} 
Assume that $\b_1 = 0$, and let $P$ be the ortogonal projection from $\L^2 (\Omega)$ onto $\X$. 
The Stokes operator defined as $A = -P \Delta$
with domain
\[ 
D(A) 
= 
\{
\u \in \J_0 \cap \H^2 (\Omega):\, \u|_{S} = 0, \,  \u\cdot \n |_{\Gamma} = 0 , \,
 [ D(\u)\n - \n \cdot ( D(\u)\n)\n]|_{\Gamma}=0
\}
\]

\noindent
 is a selfadjoint, positive definite operator
 with compact inverse.
 \end{lemma}

Observe that, since $\b_1 = 0$, the expressions for the boundary conditions become linear since the viscosity term drops out.

Moreover, from definition of $A$, it follows that $Dom(A^{1/2})=\J_0$ and $\|\ \cdot \|_{\H}$ and $|A^{1/2} (\cdot) |$ are equivalent norms in $\J_0$.
For $0 \leq \gamma_1 \leq \gamma_2 \leq 1$, we have that $ Dom(A^{\gamma_2} )\subset  Dom(A^{\gamma_1})$ and
$|A^{\gamma_1} \v | \leq C_{\gamma_1, \gamma_2} |A^{\gamma_2} \v |$, 
with $C_{\gamma_1, \gamma_2}$
independent of $\v \in Dom(A^{\gamma_2})$.
 
\vspace{0.1cm}
 An immediate consequence of this last lemma, see Brezis \cite{brezis},  is the following:
 \begin{lemma} 
\label{StokesOperatorEigenvaluesEigenvectors}
 The Stokes operator $A$ of Lemma \ref{Lema 1.8}  has a sequence
eigenvalues $\{\alpha_i\}_{i=1}^\infty$ of eigenvalues satisfying $0<\alpha_1\leq\alpha_2\leq ...$ and $\lim_{i\rightarrow +\infty}\alpha_i=+\infty$, whose associated eigenfunctions $\{\mathbf{w}^i\}_{i=1}^\infty$ form a complete orthonormal system in $\X$, $\J_0$ and $D(A)$. 
\end{lemma}

We observe that each eigenfunction $\mathbf{w}^i$ satisfies the following boundary condition
\begin{equation}
\label{bc0EigenvectorsVelocity}
\left\{\begin{array}{l}
\mathbf{w}^i = 0 \qquad {\rm on} \qquad (0,T] \times S, \\
\mathbf{w}^i \cdot \n = 0 \qquad {\rm on} \qquad (0, T] \times \Gamma,\\
\nu(m)[D(\mathbf{w}^i)\n - \n\cdot(D(\mathbf{w}^i)\n)\n] =\b_1
\quad \mbox{on} \quad (0, T] \times \Gamma, 
\end{array}\right.
\end{equation}

We will use the orthogonal projections
$P_n : \X \rightarrow span \langle \w^1, \w^2, \ldots, \w^n\rangle$, for each positive integer $n$.
When $\v \in \X$, we have $\v = \sum_{i=1}^\infty(\v, \w^i)\w^i$
and $P_n \v=\sum_{i=1}^n (\v,\w^i)\w^i$.

For $\v \in Dom(A^\gamma)$, $ 0 \leq \gamma \leq 1$, we then have that
 $A^\gamma \v =  \sum_{i=1}^\infty( \v ,\w^i) \alpha_i^\gamma \w^i$,
and denoting $\v^n = P_n \v= \sum_{i=1}^n (\v,\w^i)\w^i$,
due to the orthonormality  of the eigenfunctions, we have that
$|A^\gamma \v^n| \leq |A^\gamma \v|$
and $A^\gamma \v^n \rightarrow A^\gamma \v$ in the norm $\| \cdot \|_{A^\gamma} = |A^\gamma (\cdot) |$.
In particular, by taking $\gamma = 0$ and  $\v \in Dom(A^0) = Dom(I) = \X$, we have $\v^n \rightarrow \v$ in $\X$.
Also, by taking $\gamma = 1/2$ and  $\v \in Dom(A^{1/2})= \J_0$, we have $\v^n \rightarrow \v$ in $\J_0$ (and so also in the norm of $\H^1 (\Omega)$;
this last convergence implies that $D (\v^n) \rightarrow D(\v)$ in the $\L^2(\Omega)$-norm.

\vspace{0.2cm}
Next, since the convection-diffusion equation for $m$ implies that  $\displaystyle\frac{d}{dt}\jnt_\Om m(x,t)= 0$,
to deal with the microorganism concentration, we will introduce suitable changes of variables; 
the following functional spaces will be useful for the resulting changed variables.

Let $Y$ de the closed subspace of $L^2 (\Omega)$ consisting of
orthogonal functions to the constants, that is
\begin{equation}
\label{DefinitionY,B}
Y=\{f \in L^2 (\Omega): \int_\Omega f(x)dx =0 \}
\quad
\mbox{ and }
\quad
B=H^1 (\Omega) \cap Y.
\end{equation}

\vspace{0.2cm}
The computations involving the microorganisms concentration will require the following technical results.
\begin{lemma} 
\label{Lema 1.4.} 
(Morrey, \cite{morrey}) Under the hypotheses of Lemma
\ref{Lema 1.2}, there exists a positive constant, $C_P$,
such that
\[
|\phi| \leq C_P |\nabla\phi|, \qquad \forall \phi \in H^1
(\Omega)\cap Y.
\]
\end{lemma}

\begin{lemma}
\label{Lema 1.9} 
(Kan-On {\it et al.} \cite{kan}, Appendix, p. 150.)
Let $\overline P$ be the orthogonal projection from $L^2(\Omega)$ onto $Y$. 
The operator $A_1 : Y\mapsto Y$, defined as $\overline P(-\theta \Delta)$
defined for $\varphi \in Y \cap H^2 (\Omega)$ satisfying $\displaystyle \theta \frac{\partial\varphi}{\partial \textbf{n}} - Un_3 \varphi\vert_{\partial\Omega} = 0$
 is a selfadjoint, positive definite operator with  compact inverse. 
\end{lemma}

Moreover, from definition of $A_1$, it follows that $Dom(A_1^{1/2})=B$ and 
\[(\theta-2\,U\,C_p)^{1/2}| \nabla \varphi |
\leq 
| A_1^{1/2}\varphi |
\leq
(\theta+2\,U\,C_p)^{1/2}| \nabla \varphi |  \quad \forall \varphi\in B .
\]

\noindent
Thus, in order to $| A_1^{1/2}\varphi |$ and $| \nabla \varphi |$ become equivalent norms in $B$, throughout this work we assume the following smallness condition on $U$:
\begin{equation}
\label{SmallnessConditionOnUAuxiliaryProblem}
2 U C_P < \theta ,
\end{equation}

\noindent
where $C_P$ is the positive constant appearing in Lemma \ref{Lema 1.4.}.

\vspace{0.2cm}
Similarly as in the case of operator $A$, for $0 \leq \gamma_1 \leq \gamma_2 \leq 1$, we have that $ Dom(A_1^{\gamma_2} )\subset  Dom(A_1^{\gamma_1})$ and
$|A_1^{\gamma_1} \varphi | \leq C_{\gamma_1, \gamma_2} |A_1^{\gamma_2} \varphi |$, 
with $C_{\gamma_1, \gamma_2}$
independent of $\varphi \in Dom(A_1^{\gamma_2})$.

Moreover, as in the case of the Stokes operator, see Giga, Miyakawa \cite{GigaMiyakawa1985}. Lemma 1.3 and Lemma 1.4, 
for any $\gamma$, $0 \leq \gamma \leq 1$, the domain $Dom (A_1^\gamma)$ is the complex interpolation space $[Y, D(A_1)]_\gamma$.
As consequence, $Dom(A_1^\gamma) = Y \cap Dom( (- \theta \Delta)^\gamma)$,
and, for any $\gamma \geq 0$, the domain $Dom(A_1^\gamma)$ is continuously embedded in $Y \cap H^{2 \gamma}(\Omega)$,
where $H^{2 \gamma}(\Omega)$ is the functional space of Bessel potentials of order $2 \gamma$ based on $L^2 (\Omega)$.
This last result implies that there is a positive constant $C_{\gamma, \Omega} $ independent of $\varphi \in Dom(A_1^\gamma)$ such that
\begin{equation}
\| \varphi \|_{H^{2 \gamma}(\Omega)} \leq C_{\gamma, \Omega} | A_1^\gamma \varphi |, \quad \forall \varphi \in Dom(A_1^\gamma)
\end{equation}

\noindent
Since Sobolev embeddings hold for Bessel potentials, we obtain in particular that 
$|\nabla \varphi |_4 \leq C  \| \varphi \| _{W^{1,4} (\Omega)}  \leq C \| \varphi \| _{H^{3/2} (\Omega)}\leq C |A_1^{3/4} \varphi |  $
 for all $\varphi \in Dom(A_1^{3/4})$.
Therefore,
\begin{equation}
\label{EstimateImportant}
\begin{array}{l}
|\nabla \varphi |_4 \leq C |A_1^{1/2} \varphi | \leq \hat{C}  |A_1 \varphi |  \quad  \forall \varphi \in Dom(A_1),
\vspace{0.2cm}
\\
|\nabla \varphi |_4 \leq C |A_1^{3/4} \varphi | \leq \hat{C}  |A_1 \varphi |  \quad  \forall \varphi \in Dom(A_1).
\end{array}
\end{equation}

\noindent
These previous estimates will be important for the computations of this work.

\vspace{0.2cm}
Again, an immediate consequence of this last lemma, see Brezis \cite{brezis},  is the following
\begin{lemma}
\label{LaplacianEigenvaluesEigenvectors}
The operator $A_1$ of Lemma \ref{Lema 1.9} has a sequence of eigenvalues, $\{\beta_i\}_{i=1}^\infty$, satisfying $0<\beta_1\leq\beta_2\leq ...$ with $\lim_{i\rightarrow +\infty}\beta_i=+\infty$, whose associated normalized eigenfunctions $\{\phi^i \}_{i=1}^\infty$ form a complete orthogonal system in $Y$, $B$ and $D(A_1)$, with their natural inner products. 
\end{lemma}

We observe that each eigenfunction $\phi^i$ satisfies the following boundary condition:
\begin{equation}
\label{bc0EigenvectorsConcentration}
\left\{\begin{array}{l}
\theta \displaystyle\frac{\partial \phi^i}{\partial \n} - U \phi^i n_3 = 0
\qquad {\rm on} \qquad  (0, T) \times \partial\Omega. 
\end{array}\right.
\end{equation}

We also observe that the eigenvalues and eigenfunctions of the previous operator $A_1$ are basically the same as the eigenvalues e eigenfunctions of
operator $-\theta \Delta : L^2 (\Omega) \mapsto L^2(\Omega)$,
with domain  
$\displaystyle D(-\theta \Delta) = \{\varphi \in L^2(\Omega) \cap H^2;\, \theta \frac{\partial\varphi}{\partial \textbf{n}} - Un_3 \varphi\vert_{\partial\Omega} = 0 \}$.
In fact, being $\{\beta_i\}_{i=0}^\infty$ and $\{ \phi^i\}_{i=0}^\infty$, respectively the eigenvalues and associated eigenfunctions of $-\theta \Delta $,
we have $\beta_0 = 0$ and $\phi_0 = 1/|\Omega|^{1/2}$.
Moreover, since $-\theta \Delta $ is self-adjoint, we have that, for $i \geq 1$, $\phi^i$ is orthogonal in $L^2 (\Omega)$ to $\phi^0$; this implies that 
$\displaystyle \int_\Omega \phi^i \, dx = 0$ for $i \geq 1$; 
that is, $\phi^i$ has mean value zero for $i \geq 1$ and so $\phi^i \in Y$ and $\overline{P} \phi^i = \phi^i$  for $i \geq 1$.
Since $-\theta \Delta \phi^i = \beta^i \phi^i$, we obtain that $\displaystyle \int_\Omega -\theta \Delta \phi^i \, dx = \beta^i \int_\Omega \phi^i \, dx = 0$ for  $i \geq 1$,
and so  $\overline{P} (-\theta \Delta \phi^i ) = -\theta \Delta \phi^i $ for  $i \geq 1$,
We conclude that $A_ 1 \phi^i = \overline{P} (-\theta \Delta ) \phi^i = -\theta \Delta \phi^i  = \beta^i \phi^i$;
that is, $\{\beta_i\}_{i=1}^\infty$ and $\{ \phi^i\}_{i=1}^\infty$, are respectively the eigenvalues and associated eigenfunctions of $A_1 $.

We will use the orthogonal projections
$\overline P_n : Y
\rightarrow span \langle \phi^1, \phi^2, \ldots, \phi^n\rangle$, for each positive integer $n$.
When $\phi \in Y$, we have $\phi = \sum_{i=1}^\infty(\phi, \phi^i)\phi^i$
and $\overline P_n(\phi)=\sum_{i=1}^n(\phi,\phi^i)\phi^i$.
For $\phi \in Dom (A^\gamma)$, $0 \leq \gamma \leq 1$, we then have that
 $A_1^\gamma \phi =  \sum_{i=1}^\infty( \phi ,\phi^i) \beta_i^\gamma \phi^i$,
and denoting $\phi^n = \overline P_n(\phi)= = \sum_{i=1}^n (\phi,\phi^i)\phi^i$,
due to the orthonormality  of the eigenfunctions, 
we have that
$|A_1^\gamma \phi^n| \leq |A_1^\gamma \phi|$
and $A_1^\gamma \phi^n \rightarrow A_1^\gamma \phi$ in the norm $\| \cdot \|_{A_1^\gamma} = |A_1^\gamma (\cdot) |$.
In particular, by taking $\gamma = 0$ and  $\phi\in Dom(A_1^0) = Dom(I) = Y$, we have $\phi^n \rightarrow \phi$ in $Y$.
Also, by taking $\gamma = 1/2$ and  $\phi \in Dom(A_1^{1/2}) \subset Y \cap H^1(\Omega)$, with continuous embedding,
we have $\phi^n \rightarrow \phi$ in the norm $\| \cdot \|_{A_1^{1/2}}$, and so also in the norm of $H^1 (\Omega)$.
By taking $\gamma = 1$ and  $\phi \in Dom(A_1)\subset Y \cap H^2 (\Omega)$, with continuous embedding, 
we have $\phi^n \rightarrow \phi$ in the norm $\| \cdot \|_{A_1}$, and so also in the norm of $H^2 (\Omega)$;
this last convergence implies in particular that $\nabla \phi^n \rightarrow \nabla \phi$ in the $L^4(\Omega)$-norm.

\vspace{0.2cm}
We will also need the following results.

\begin{lemma} 
\label{helmholtz}
 Let $\u\in \J_0\cap \H^2 $ and consider the Helmholtz decomposition of $-\Delta \u$, i.e. 
\[
-\Delta \u=A\u+\nabla q,
\]

\noindent
where $q\in H^1 (\Omega)$ is taken such that $\displaystyle\int_\Omega q\; dx=0$.
 Then, there exists a positive constant $C$ and, for any $\eps>0$, an associated positive constant $C_\eps$, independent of $\u$ , such that
\[
\| q \|_1\leq C | A\u | \quad \mbox {and }\quad
| q |\leq C_\eps | \nabla \u |+\eps | A\u |, \qquad \forall \u \in D(A).
\]
\end{lemma}

The proof of this previous lemma is similar to the one of Lemma 3.4 presented in Lorca, Boldrini \cite{lorca2}.

\begin{lemma}
 \label{blanca} 
Let  $a(t)\geq 0$ be a non-decreasing function
 and $b(t)\geq 0$ be a summable function in $[0,T]$. 
Assume that the integral inequality
\[
\zeta(t) + \int_0^t \zeta^*(s)ds \leq a(t) +\int_0^t b(s)\zeta(s)ds
\]

\noindent
holds for nonnegative continuous functions $\zeta$ and $\zeta^*$ on $[0,T]$.
Then
\[
\displaystyle
\zeta(t)+\int_0^t\zeta^*(s)ds \leq a(t)  e^{\int_0^t b(s)ds}.
\]
\end{lemma}

This previous lemma is a variant of Gronwall inequality; see for instance Hale \cite{hale} .

\begin{lemma} 
\label{lema110}
 Let $\varphi, \psi $ positive functions such that
\[
\left\{\begin{array}{l}
\varphi'(t)+\psi(t)\leq C\varphi^3(t)+h(t)\varphi(t)+g(t), \quad t\in[0,T]\\
\varphi(0)=\varphi_0,
\end{array}
\right.
\]
where $\varphi_0\geq 0$ and, $h(t)$ and $g(t)$ are positive and continuous functions. Let $T_1$, $0\leq T_1\leq T$, such that 
\[
C\int_0^{T_1}(2\varphi_0+2\int_0^{s}g(\tau)\,d\tau)^2+\int_0^{T_1}h(s)\,ds\leq \frac{1}{2}\ln{\frac{3}{2}}. 
\]
Then, 
\[
\varphi(t)+\int_0^{t}\psi(t)\,ds\leq 2(\varphi_0+\int_0^{t}g(s)\,ds)
\]
for all $t\in[0,T_1]$.
\end{lemma}

This previous result is consequence of Gronwall lemma and appears in Lorca, Boldrini \cite{lorca2}.

\section{Existence of global weak evolution solutions}
\label{ExistenceOfGlobalEvolutionaryWeakSolutions}

By doing standard computations, using Lemma \ref{lema7}, we obtain the following weak formulation:
\begin{definition} 
\label{DefinitionEvolutiveWeakSolution}
A pair of functions $(\u, m)$ is a weak solution of \eqref{sistema2}-\eqref{ic}-\eqref{bc0} if 
\[
\begin{array}{ll}
\u \in L^2(0,T;\J_0) \cap L^\infty (0,T; \X), & 
m \in L^2 (0,T; H^1 (\Omega) \cap L^\infty (0,T; L^2 (\Omega)),
\vspace{0.1cm}
\\
 \partial_t \u \in  L^1 (0, T; \J_0'), &
 \partial_t m \in L^1 (0, T;  B') ,
 \end{array}
\]

\noindent
where $B$ is the functional space defined in (\ref{DefinitionY,B});
for all $\w \in \J_0$ and $\phi \in B$, the following identities hold in $\mathcal{D}'(0,T)$:
\begin{equation}
\label{WeakFormulation}
  \left\{
\begin{array}{l}
\displaystyle
\left< \partial_t \u, \w  \right> 
 + 2(\nu(m)D(\u), D(\w)) 
+ (\u \cdot\nabla \u,\w) 
 =  
 - (m \cdot \chi, \w)
+ (\f , \w)
+ 2 ( \b_1,  \w)_\Gamma
  \\
  \displaystyle
 \left< \partial_t m, \phi \right> + \theta (\nabla m, \nabla \phi) 
 + (\u \cdot\nabla m,\phi) 
 - U\left( m , \frac{\partial\phi}{\partial x_3} \right) 
 = 
 0,
\\
\u = \u_0  \quad \mbox{and} \quad m = m_0, \quad \mbox{at time} \; \; t=0 .
\end{array}
\right.
\end{equation}

\noindent
In \eqref{WeakFormulation}, $\left< \cdot, \cdot \right>$ denotes the duality action between a space and its dual.
\end{definition}

For this problem, we have the following result

\begin{theorem} 
\label{WeakSolutionExistence} 
Assume that $\f \in L^2 (0,T; \X)$,
$\b_1 \in L^2 (0,T; L^2(\Gamma) )$ 
and that the initial conditions satisfy $\u_0 \in \X $ and  $m_0 \in  L^2(\Omega)$.
Assume also that  $\nu (\cdot)$ is a continuous function  satisfying	
\begin{equation}
\label{nus}
\nu_0=\inf\{\nu(m), m\in \mathbb{R}\}>0,\qquad
\nu_1=\sup\{\nu(m), m\in \mathbb{R}\}<+\infty.
\end{equation} 

 \noindent
Then, there exists a weak solution, $(\u, m)$, of  \eqref{sistema2}-\eqref{ic}-\eqref{bc0} according to Definition \ref{DefinitionEvolutiveWeakSolution}.
 Moreover,
\begin{equation}
\label{limci}
|\u(t)-\u_0| \rightarrow 0 \quad {\rm and} \quad |m(t)-m_0|
\rightarrow 0 \quad \mbox{as} \; t \rightarrow 0^+.
\end{equation}
\end{theorem}

\subsection{Proof Theorem  Theorem \ref{WeakSolutionExistence}}
To prove the existence of global weak evolution solutions, it is convenient to introduce the following change of variables
\begin{equation}
\label{ChangeOFVariablesConcentration}
\zeta= m- \frac{\alpha}{|\Omega|},
\end{equation}

\noindent
which rewrites \eqref{TotalMassConservation}
\[
\int_\Omega \zeta (x, t) dx = 0, \; \mbox{for all} \; t \in [0,T] .
\]

Then, in terms of  $\u$ and $\zeta$, problem (\ref{WeakFormulation}) is rewritten as follows: 
find  $(\u, \eta)$ satisfying
\begin{equation}
\label{RegularityOfWeakSoltuinsProposition}
\begin{array}{ll}
\u \in L^2(0,T;\J_0) \cap L^\infty (0,T; \X), & 
\zeta \in L^2 (0,T; B) \cap L^\infty (0,T; Y),
\vspace{0.1cm}
\\
 \partial_t \u \in  L^1 (0, T; \J_0'), &
 \partial_t \zeta \in L^1 (0, T;  B'),
 \end{array}
\end{equation}

\noindent
and such that, for all $\w \in \J_0$ and $\phi \in B$, the following identities hold:
\begin{equation}
\label{WeakFormulationOtherProposition}
  \left\{
\begin{array}{l}
\displaystyle
 \left< \partial_t \u, \w \right> 
 + 2\left(\nu(\zeta + \frac{\alpha}{|\Omega|})D(\u), D(\w)\right) 
+ (\u \cdot\nabla \u,\w) 
 + ( \zeta  \chi, \w)
\\
\hspace{1.5cm}
\displaystyle
 =  
  - \frac{\alpha}{|\Omega|}   ( \chi, \w)
+ (\f , \w)
+ 2 ( \b_1,  \w)_\Gamma ,
  \\
  \displaystyle
 \left< \partial_t \zeta, \phi \right> + \theta (\nabla \zeta, \nabla \phi) 
 + (\u \cdot\nabla \zeta,\phi) 
 - U\left( \zeta  , \frac{\partial\phi}{\partial x_3} \right) 
 =
   U \frac{\alpha}{|\Omega|} \left( 1 , \frac{\partial\phi}{\partial x_3} \right) 
,
\end{array}
\right.
\end{equation}

\noindent
in $\mathcal{D}'(0,T)$, and also
\begin{equation}
\label{InitialConditionOtherProposition}
\begin{array}{l}
\displaystyle
\u(t) \rightarrow \u_0 \quad \mbox{strongly in} \; \X \; {\rm as} \; t \rightarrow 0^+ ,
\\
\displaystyle
\zeta(t) \rightarrow \zeta_0 = m_0 - \frac{\alpha}{|\Omega|} \quad \mbox{strongly in} \; Y \; {\rm as} \; t \rightarrow 0^+. 
\end{array}
\end{equation}

Obviously, once we have a solution of this problem, we obtain a solution of (\ref{WeakFormulation}).
Thus, it is enough to prove the following proposition:
\begin{proposition} 
\label{Proposition 3.1.} 
Under the conditions of Theorem \ref{ssolestac} and  $(\v_0, \eta_0) \in X \times Y$,
 there exists a weak solution $(\u, \zeta)$ of  \eqref{WeakFormulationOtherProposition} and \eqref{InitialConditionOtherProposition},
 with the regularities  stated in \eqref{RegularityOfWeakSoltuinsProposition}. Moreover, \eqref{InitialConditionOtherProposition} is also true.

\end{proposition}

\noindent
{\bf Proof:}
The following proof is similar to the corresponding one in Kan-On {\it et al.} \cite{kan}.

We fix the Schauder bases
$(\overline \w^j)_1^\infty$ and $(\overline
\phi^j)_1^\infty$, of  $\J_0$ and $B$, respectively, and build a weak solution of   \eqref{WeakFormulationOtherProposition} and (\ref{InitialConditionOtherProposition}) by using the Galerkin approximations:
\begin{equation}
\label{galerkinsol}
\u^n(t,x) = \sum^n_{j=1} \overline{c}_{n,j} (t) \overline \w^j(x);
\quad \zeta^n (t,x) = \sum^n_{\ell=1} \overline{d}_{n,\ell} (t)
\overline \phi^\ell(x),
\end{equation}

\noindent
satisfying the following approximate problem:
\begin{equation}
\label{ApproximateEquations}
 \left\{
\begin{array}{l}
\displaystyle
(\partial_t \u^n, \overline \w^j)
  + 2\left(\nu(\zeta^n +  \frac{\alpha}{|\Omega|} )D(\u^n), D(\overline \w^j)\right)
   + (\u^n\cdot\nabla \u^n, \overline \w^j)
 + (\zeta^n \chi, \overline \w^j)
 \\
 \hspace{2cm}
 \displaystyle
 =   - \frac{\alpha}{|\Omega|}   ( \chi, \overline{\w}^j )
+ (\f , \overline{\w}^j)
+ 2 ( \b_1,  \overline{\w}^j)_\Gamma ,
\vspace{0.1cm}
  \\
  \displaystyle
( \partial_t \zeta^n, \overline
\phi^\ell)
  + \theta (\nabla \zeta^n, \nabla \overline \phi^\ell)
  + (\u^n\cdot\nabla \zeta^n, \overline \phi^\ell)
- U\left(\zeta^n, \frac{\partial \overline{\phi}^\ell }{\partial x_3} \right) 
=    U \frac{\alpha}{|\Omega|} \left( 1 , \frac{\partial \overline{\phi}^j }{\partial x_3} \right)  ,   
 \vspace{0.1cm}
 \\
 \u^n =  \u_0^n , \quad \ \zeta^n = \zeta_0^n \quad \mbox{at time} \; t = 0,
\end{array}
\right.
\end{equation}

\noindent
for all $1 \leq j, \ell \leq n$. 

Here, for each $n$, $ \u_0^n  \in span \{\overline{\w}^1, \ldots, \overline{\w}^n \}$ and  $ \zeta_0^n  \in span \{\ \overline{\phi}^1, \ldots, \overline{\phi}^n \}$ such that $ \u_0^n \rightarrow \u_0$ strongly in $X$ and $ \zeta_0^n \rightarrow \zeta_0$ strongly in $Y$.
In particular, one could take $ \u_0^n = Q_n \u_0$ and  $ \zeta_0^n = \overline Q_n \zeta_0$,
with $Q_n : X \rightarrow span \{\overline{\w}^1, \ldots, \overline{\w}^n \}$ and $\overline Q_n Y : \rightarrow span \{\ \overline{\phi}^1, \ldots, \overline{\phi}^n \}$ orthogonal projections.

It is standard that this last system has a maximal solution defined in some
interval $[0, t_n)$. The following estimates show that actually, $t_n = T$.

By multiplying the first and second equations in (\ref{ApproximateEquations}) respectively by $\overline{c}_{n,j}$ and $\overline{d}_{n,l}$, 
adding with respect to  $j$ and  $l$, we obtain:
\begin{eqnarray}
&& \frac{1}{2} \frac{d}{dt} |\u^n|^2 
+ 2 \left(\nu(\zeta^n +  \frac{\alpha}{|\Omega|}) D(\u^n), D(\u^n)\right) 
+ (\zeta^n  \chi, \u^n) 
\nonumber
\\
&&
 \hspace{2cm}
 \displaystyle
 =   
 - \frac{\alpha}{|\Omega|}   ( \chi, \u^n )
+ (\f , \u^n)
+ 2 ( \b_1,  \u^n)_\Gamma ,
\label{IdentityU1}
\vspace{0.2cm}
\\
&&\frac12 \frac{d}{dt} |\zeta^n|^2 
+ \theta |\nabla\zeta^n|^2 
 - U\left(\zeta^n, \frac{\partial \zeta^n}{\partial x_3}\right) 
 = 
     U \frac{\alpha}{|\Omega|} \left( 1 , \frac{\partial \zeta^n }{\partial x_3} \right) .
 \label{IdentityEta1}
\end{eqnarray}

By using (\ref{nus}), H\"older and Young inequalities and standard Sobolev embedding's we have the following estimates: 
\begin{eqnarray*}
\frac12 \frac{d}{dt} |\u^n|^2 
+ 2 \nu_0 |D(\u^n)|^2 & \leq & 
 |\zeta^n| \ |\u^n|  
 + \frac{\alpha}{|\Omega|^{1/2}}   |\u^n |
+ |\f |  |\u^n|
+ 2 \| \b_1\|_{L^2 (\Gamma ) }   \| \u^n \|_{L^2 (\Gamma) } , 
\\
& \leq &\frac{3}{2}\nu_0 |D(\u^n)|^2
+ C_{\nu_0}  | \nabla \zeta^n|^2 
+ C_{\nu_0} \left(  |\f |^2+  \frac{\alpha^2}{|\Omega| } + \| \b_1 \|_{L^2 (\Gamma) }^2 \right)
,
\\
\frac12 \frac{d}{dt} |\zeta^n|^2 + \theta |\nabla\zeta^n|^2 & \leq &
U|\zeta^n| \ |\nabla\zeta^n|
+ U \frac{\alpha}{|\Omega|^{1/2}} | \nabla \zeta^n | .
\\
& \leq &
\frac{\theta}{2}  |\nabla\zeta^n|^2
+C_{\theta U}  |\zeta^n|^2 + C_{\theta U} \frac{\alpha^2}{|\Omega|}  .
\end{eqnarray*}

These estimates imply that
\begin{eqnarray}
\frac{d}{dt} |\u^n|^2 +  \nu_0 |D(\u^n)|^2 & \leq & 
2 C_{\nu_0}  | \nabla \zeta^n|^2 
+2 C_{\nu_0} \left(  |\f |^2+  \frac{\alpha^2}{|\Omega| } + \| \b_1 \|_{L^2 (\Gamma) }^2 \right) ,
\label{FirstInequalityApprox}
\\
 \frac{d}{dt} |\zeta^n|^2 + \theta |\nabla\zeta^n|^2 & \leq & 2C_{\theta U}  |\zeta^n|^2 + 2C_{\theta U} \frac{\alpha^2}{|\Omega|} .
 \label{SecondInequalityApprox}
\end{eqnarray}

By multiplying the second of these last estimates by a positive constant $C_0$ such  that
\[
\theta C_0 = 4 C_{\nu_0}    ,
\]

\noindent
and then adding the resulting inequality to the first estimate, we obtain:

\[
\begin{array}{c}
\fra{d}{dt} (|\u^n|^2 + 4 C_{\nu_0} |  \zeta^n|^2) + \nu_0 |D(\u^n)|^2 +  2 C_{\nu_0}  |\nabla\zeta^n|^2 
\\
\displaystyle
\leq 
2 C_{\nu_0} \left(  |\f |^2+  \frac{\alpha^2}{|\Omega| } + \| \b_1 \|_{L^2 (\Gamma) }^2 \right)
+  8 C_{\nu_0}   C_{\theta U} (   |\zeta^n|^2 +  \frac{\alpha^2}{|\Omega|}  )
\\
\displaystyle
\leq 
  C_1 ( |\u^n|^2 + 4 C_{\nu_0}   |\zeta^n|^2 )
+ C_2 \left(  |\f |^2+  \frac{\alpha^2}{|\Omega| } + \| \b_1 \|_{L^2 (\Gamma) }^2 \right) .

\end{array}\]

By integrating this last inequality in $[0,t]$, applying Gronwall's inequality and observing that $\|Q_n\| \leq 1, \ \|\overline{Q}_n\| \leq 1$, it
follows that
\begin{equation}
\label{G}
 |\u^n(t)|^2 + 4 C_{\nu_0}   |\zeta^n(t)|^2 
+ \nu_0 \int^t_0 |D(\v^n) (s)|^2 ds 
+  2 C_{\nu_0}  \int^t_0 |\nabla \zeta^n (s)|^2 ds 
\leq 
G(t),
\end{equation}

\noindent
where 
\begin{equation}
\label{G(t)}
G(t) = \left( |\u_0|^2 + 4 C_{\nu_0}   |\zeta_0|^2  + \int_0^t C_2 \left(  |\f |^2+  \frac{\alpha^2}{|\Omega| } + \| \b_1 \|_{L^2 (\Gamma) }^2 \right) ds \right) \exp (C_1 t) ,
\end{equation}

\noindent
which is a continuous function, independent of $n$ and bounded on bounded intervals.

Thus, the sequence $(\u^n, \zeta^n)$ is bounded in 
$L^2 (0, T; \J_0) \cap L^\infty (0, T; X) \times L^2  (0, T;B) \cap L^\infty (0,T;Y)$.

To obtain  estimates to $\partial_t \u^n$ and $\partial_t \zeta^n$,
we use arguments analogous to those of Lions \cite{lions1}.
For this, we define the following operators:
\[
\begin{array}{l}
\displaystyle
 \langle A_\nu(u, \zeta), \w\rangle = 2 \left(\nu(\zeta^n + \frac{\alpha}{|\Omega|}) D(\u), D(\w) \right),  
\\
 \langle F(\u,\v), \w\rangle = (\u\cdot\nabla \v, \w), 
\\
\displaystyle
 \langle H(\zeta), \w\rangle = (\zeta \chi , \w)  + ( \frac{\alpha}{|\Omega|}  \chi, \w)
 - (\f , \w)
- 2 ( \b_1, \cdot \w)_\Gamma , 
\vspace{0.1cm}
\\
\displaystyle
 \langle \overline A(\zeta), \phi\rangle 
 = \theta (\nabla \zeta , \nabla\phi)  ,
 \\
 \langle \overline F(\u, \zeta),\phi\rangle = (\u\cdot\nabla \zeta, \phi),
 \\
 \displaystyle
  \langle \overline{H}(\zeta), \w\rangle =  - U \left( \zeta ,  \frac{\partial\phi}{\partial x_3}\right)  - U \left( \frac{\alpha}{|\Omega|}, \frac{\partial\phi}{\partial x_3}\right) .

\end{array}
\]

\noindent
which satisfy the following estimates:
\begin{equation}
\begin{array}{c}
||A_\nu(\u, \zeta)||_{\J_0'} \leq 2 \nu_1 |D(\u)|,
\qquad
\|F(\u, \v)\|_{\J_0'} \leq C |D(\u)| \ |D(\v)|, \label{AFH0}
\vspace{0.1cm}
\\
\displaystyle
\|H(\zeta)\|_{\J_0'} \leq C \left( |\zeta| +  \frac{\alpha}{|\Omega|}+ |\f | + 2 \|\b_1 \|_{L^2 (\Gamma)} \right),
\end{array}
\end{equation}
\begin{equation}
\begin{array}{c}
\displaystyle
\|\overline A(\zeta)\|_{B'} \leq \theta |\nabla \zeta| ,
\qquad
\|\overline F(\u, \zeta)\|_{B'} \leq C|D(u)| \ |\nabla \zeta| ,
\vspace{0.1cm}
\\
\displaystyle
\| \overline{H}(\zeta)\|_{\J_0'} \leq C  U\left( |\zeta| + \frac{\alpha}{|\Omega|} \right)
.
\label{AFH}
\end{array}
\end{equation}

\noindent
From the definition of weak approximate solution, we then have that 
\begin{eqnarray}
\partial_t \u^n & = & 
-Q^*_n 
\left(
A_\nu(\u^n , \zeta^n) 
+ F(\u^n, \u^n)
+ H(\zeta^n ) 
\right)
\qquad
{\rm in} \; \J_0',
\label{pt0}
\\
\partial_t \zeta^n & = &
 - \overline{Q}_n^*
 \left(
 \overline A(\zeta^n) 
 + \overline F(\v^n, \zeta^n)
 + \overline{H}(\zeta^n 
 \right)
 ) 
\quad \mbox{in} \; B',
\label{pt}
\end{eqnarray}

\noindent
where $Q^*_n$ and $\overline{Q}^*_n$ are respectively the adjoint operators of the orthogonal projection $Q_n$ and $\overline Q_n$.

Equation  (\ref{pt0}), estimates (\ref{AFH0}) and the previously proved estimates imply that
$\{\partial_t \u^n\}$ is a bounded sequence in $L^1 (0,T;\J_0')$, 
and so $\{\u^n\}$ is relatively compact in $L^q(0,T;\X)$, for all $1 \leq q < \infty$ (Simon \cite{simon}).

Analogously, from (\ref{pt}), (\ref{AFH}) and the previously proved estimates, 
we get that $\{\partial_t \zeta^n\}$
is bounded  in $L^1 (0, T;B')$ ,
and so $\{\zeta^n\}$ is relatively compact in $L^q(0, T; Y)$, for all $1 \leq q < \infty$.

Therefore, we conclude that there is a subsequence of
$\{(\u^n,\zeta^n)\}$, which for simplicity of exposition we do not relabel, such that
\begin{eqnarray*}
&& (\u^n, \zeta^n) \ \rightarrow \ (\u, \zeta) \quad \mbox{weakly in} \quad L^2 (0, T; \J_0 \times B), 
\\
&& (\u^n, \zeta^n) \ \rightarrow \ (\u, \zeta) \quad \mbox{$\star$-weakly in} \quad L^\infty (0, T; \X \times Y) ,
\\
&& (\u^n, \zeta^n) \ \rightarrow \ (\u, \zeta) \quad \mbox{strongly in} \quad L^2 (0, T; \X \times Y) ,
\\
&& (\partial_t \u^n, \partial_t \zeta^n) \ \rightarrow \ (\partial_t \u, \partial_t \zeta) \quad \mbox{weakly in} \quad L^1 (0, T; \J_0' \times B').
\end{eqnarray*}
Moreover, we can choose the subsequence such that
\[
\begin{array}{l}
(D(\u^n), \nabla\zeta^n) \ \rightarrow \ (D(\u), \nabla\zeta) \quad \mbox{weakly in} \: L^2 (0, T; \L^2 \times L^2) ,
\vspace{0.2cm}
\\
\zeta^n \rightarrow \zeta \quad \mbox{a.e. in} \: [0, T] \times \Omega.
\end{array}
\]

These convergences suffice to pass to limit in the equations of (\ref{ApproximateEquations}) to obtain 
\eqref{WeakFormulationOtherProposition}, and thus also (\ref{WeakFormulation}).
In fact, the only non-standard term to pass to the limit is $\displaystyle \nu(\zeta^n +  \frac{\alpha}{|\Omega|} )D(\u^n)$.
To deal with this term, we observe that, from the continuity of $\nu (\cdot)$, 
we have that
\[
\displaystyle \nu(\zeta^n (\cdot) +  \frac{\alpha}{|\Omega|} ) - \nu(\zeta (\cdot) +  \frac{\alpha}{|\Omega|} )  \rightarrow 0 \quad \mbox{a.e. in} \: [0, T] \times \Omega.
\]

\noindent
This, the fact that $\displaystyle \left|\nu(\zeta^n (\cdot) +  \frac{\alpha}{|\Omega|} ) - \nu(\zeta (\cdot) +  \frac{\alpha}{|\Omega|} ) \right|^2 \leq 4 \nu_1^2$,
and since $\Omega$ is a bounded domain, by using the dominated convergence theorem, then imply that 
\[
\displaystyle \nu(\zeta^n (\cdot) +  \frac{\alpha}{|\Omega|} )  \rightarrow  \nu(\zeta (\cdot) +  \frac{\alpha}{|\Omega|} )  
\quad \mbox{strongly in} \: L^2 (0, T; L^2 (\Omega)) .
\]

\noindent
Therefore, $\displaystyle \nu(\zeta^n +  \frac{\alpha}{|\Omega|} )D(\u^n) \rightarrow  \nu(\zeta +  \frac{\alpha}{|\Omega|} )D(\u)$ weakly in $L^2 (0, T; L^2 (\Omega)) $, and we can pass the limit in this term.

\vspace{0.1cm}
To prove (\ref{limci}), note that we have that $\u \in C([0,T]; \J_0')$ t since $\partial_t \u \in L^1 (0, T; \J_0')$;
thus, $\u (t) \rightarrow \u(0)= \u_0$ weakly in $\J_0$ as $t \rightarrow 0+$.
 
 From this result, the fact that $\u^n$ is uniformly bounded in $L^\infty (0, T; \X) $ and also that $\J_0$ is dense in $\X$,
 we also obtain that $\u(t) \rightarrow \u_0$  weakly in $ \X $ as $ t \rightarrow 0^+ $.
 In fact, being $\X$ reflexive, it is enough to show that
 $(\u(t) , \w) \rightarrow (\u_0, \w)$ for any $\w \in \X$ as $ t \rightarrow 0^+ $.
 Then, given a $\w \in \X$ and an arbitrary  $\delta > 0$, there is a $\w_\delta \in \J_0$ such that$|\w - \w_\delta |\leq \delta$.
 This implies that 
 $| (\u(t) - \u_0, \w) | 
 \leq | 
 (\u(t) - \u_0, \w_\delta) | + | (\u(t) - \u_0, \w - \w_\delta) |  
 \leq 
 | (\u(t) - \u_0, \w_\delta) | + (\|\u(t) \|_{\L^\infty \X} + |\u_0|) \delta) 
  \leq 
 | (\u(t) - \u_0, \w_\delta) | + C \delta)$.
 By taking the $\limsup$ as $ t \rightarrow 0^+ $  of both sides of this last inequaltity, using the weak convergence in $\J_0$, we then obtain
  $ \limsup_{t \rightarrow 0^+ } | (\u(t) - \u_0, \w) | \leq C \delta $.
   Since $\delta$ is arbitrary, we obtain $ \lim_{t \rightarrow 0^+ } | (\u(t) - \u_0, \w) | = 0$ as required.

 Thus, since $X$ is a uniformly convex space 
and $|\cdot|^2$ is a lower semicontinuous function in the weak
topology, we conclude that, as $t \rightarrow 0^+$,
\begin{equation}
\label{L2ContinuityTimeOrigin1}
|\u_0|^2 \leq \lim\inf |\u(t)|^2.
\end{equation}

On the other hand, by taking $\w = \u (\cdot, t)$ in the first equation of \eqref{WeakFormulationOtherProposition},
using that $\displaystyle \left< \partial_t \u, \u \right> =  \frac{1}{2} \frac{d}{dt} |\u|^2$,
see Temam \cite{temam} p. 176,
and proceeding similarly as we did to obtain \eqref{FirstInequalityApprox}, we get that
\[
\frac{d}{dt} |\u|^2 +  \nu_0 |D(\u)|^2  
\leq  
2 C_{\nu_0}  | \nabla \zeta|^2 
+2 C_{\nu_0} \left(  |\f |^2+  \frac{\alpha^2}{|\Omega| } + \| \b_1 \|_{L^2 (\Gamma) }^2 \right) .
\]

By integrating on time this last inequality, we get that
\[
|\u (t)|^2 \leq |\u_0|^2 + C \int^t_0 (2 C_{\nu_0}  | \nabla \zeta|^2 
+2 C_{\nu_0} \left(  |\f |^2+  \frac{\alpha^2}{|\Omega| } + \| \b_1 \|_{L^2 (\Gamma) }^2 \right) )  ds\quad \mbox{a.e. in } [0,T],
\]

\noindent
which implies that, as $t \rightarrow 0^+$, 
\begin{equation}
\label{L2ContinuityTimeOrigin2}
\lim\sup |\u (t)|^2 \leq |\u_0|^2.
\end{equation}

\noindent
From \eqref{L2ContinuityTimeOrigin1} and \eqref{L2ContinuityTimeOrigin2}, we obtain that
$|\u(t)| \rightarrow |\u_0|$ as $t \rightarrow 0^+$. 
By combining this last result with the weak convergence (Brezis \cite{brezis}), we finally get that
\[
\u (t) \rightarrow \u_0 \quad \mbox{strongly in} \; \X  \; \mbox{as} \; t \rightarrow 0^+.
\]

Next, by working similarly as we just did for the $\u$, taking $\phi = \eta (\cdot, t)$ in the second equation of \eqref{WeakFormulationOtherProposition},
we also obtain that
\[
\zeta(t) \rightarrow \zeta_0  = m_0 - \frac{\alpha}{|\Omega|} \quad \mbox{strongly in} \;  Y  \; \mbox{as} \; t \rightarrow 0^+. 
\]

These last two results are exactly  \eqref{InitialConditionOtherProposition}, which gives (\ref{limci}), 
and thus the proof of Proposition \ref{Proposition 3.1.}  is complete.
\hfill $\Box$

\section{Existence of local strong evolution solutions}
\label{sec:4}

Imposing further conditions, it is possible to improve the regularity of the weak solutions obtained in Theorem  \ref{WeakSolutionExistence}, at least for small time intervals, and also to prove their uniqueness. In fact, we have the following result:
\begin{theorem} 
\label{localsolMainResults}
Assume that  $\f \in L^2 (0,T; \X)$, $\b_1 = 0$ and that the initial conditions satisfy $\u_0 \in \J_0$ and $m_0 \in H^2 (\Omega)$. 
Assume also that $\nu (\cdot)$ is a function of class  $C^1$ satisfying (\ref{nus}) and
\[
\displaystyle\sup_{\mathbb{R}} |\nu'(m)| < \nu_1' < \infty.
\]

\noindent
Then, for $U$ sufficiently small, there exists $T^* = T^* (\Omega, \u_0, m_0, \nu) \leq T$,
such that problem  (\ref{sistema3})-(\ref{ic3}) has a unique solution $(\u, m)$
satisfying
\begin{eqnarray}
&& \u \in L^\infty (0, T^*, \J_0)   \cap L^2 (0, T^*; \J_0 \cap \H^2 (\Omega)),
\quad
\displaystyle
\partial_t \u \in L^2 (0, T^*; \L^2 (\Omega)),
\label{r1}
\\
&& m \in L^\infty (0, T^*; H^2 (\Omega)) ,
\quad 
\displaystyle
\partial_t m \in L^\infty (0, T^*; L^2 (\Omega)) ;
\label{r2} 
\end{eqnarray}

\noindent
Moreover,
\begin{equation}
\label{142}
\begin{array}{l}
 \u(t) \rightarrow \u_0 \ \mbox{strongly in} \ \J_0  \; \mbox{as} \; t \rightarrow 0^+,
\vspace{0.2cm}
\\
m(t) \rightarrow m_0 \ \mbox{strongly in} \ H^1 (\Omega) \, \mbox{and weakly in} \ H^2 (\Omega) \; \mbox{as} \; t \rightarrow 0^+,
\end{array}
\end{equation}

\noindent
and
\begin{equation} 
\left\{\begin{array}{l} 
\displaystyle
P \left(
\partial_t \u
- 2 \ {\rm div}\, (\nu(m)D(\u))
+ \u\cdot \nabla \u 
+m \chi - \f 
\right) 
= 0 \quad {\rm in} \ L^2(0,T^*;\X),
\vspace{0.1cm}
\\
\displaystyle
\partial_t m - \theta \Delta m 
+ \u\cdot \nabla m 
+ U \displaystyle\frac{\partial m}{\partial x_3} 
=
0 
\quad {\rm in} \ L^\infty(0, T^*;L^2 (\Omega)).
\end{array}\right.
\end{equation}

\end{theorem}

\subsection{Proof of Theorem \ref{localsolMainResults}}

We starting by stressing that from now on we assume that $\b_1 = 0$ and introducing the following auxiliary problem.

Given a fixed $\alpha > 0$, let us consider the solution $m_\alpha (\x)$ of the following system:
\begin{equation}
\label{AuxilaryProblemMalpha}
\left\{
\begin{array}{l}
\displaystyle
 - \theta \Delta m_\alpha + U \frac{\partial m_\alpha}{\partial x_3} = 0 \qquad {\rm in} \; \Omega ,
\vspace{0.2cm}
\\
\displaystyle
\theta \displaystyle\frac{\partial m_\alpha}{\partial \n} - Um_\alpha n_3 = 0
\qquad {\rm on} \;  \partial\Omega ,
\vspace{0.2cm}
\\
\int_\Omega m_\alpha \, dx = \alpha .
\end{array}
\right.
\end{equation}

For problem (\ref{AuxilaryProblemMalpha}), we have the following result:
\begin{lemma}
\label{AuxiliaryProblemLemma}
Assume that  $U$ is small enough such that (\ref{SmallnessConditionOnUAuxiliaryProblem}) holds.
Then, there exists a unique solution $m_\alpha$ of the auxiliary problem (\ref{AuxilaryProblemMalpha}).
Moreover, 
\begin{equation}
\label{MalphaEstimates}
\begin{array}{l}
\displaystyle
|\nabla m_\alpha | \leq U \frac{2 \alpha}{\theta |\Omega|^{1/2}}, 
\qquad 
 \| m_\alpha \|_1  \leq  C \, \left[ \, \alpha + U \frac{2 \alpha}{\theta |\Omega|^{1/2}} \, \right],
\vspace{0.2cm}
\\
\displaystyle
 | \Delta m_\alpha |  \leq \frac{2 U^2  \alpha}{\theta^2 |\Omega|^{1/2}},
\qquad \;
\| m_\alpha \|_ 2 \leq C [ \,  \, \left[  \alpha  + \left( 1+   \frac{U}{\theta} \right) U \frac{2 \alpha}{\theta |\Omega|^{1/2}} \, \right]  ,
\vspace{0.2cm}
\\
\displaystyle
|\nabla m_\alpha |_4 \quad \mbox{decreases to zero as $ U$ approaches zero}.
\end{array}
\end{equation}

\end{lemma}

\noindent
{\bf Proof:}
By using the following change of variables,
\begin{equation}
\label{ChangeOFVariablesConcentrationAuxiliary}
\eta_\alpha = m_\alpha - \frac{\alpha}{|\Omega|},
\end{equation} 

\noindent
which is similar to \eqref{ChangeOFVariablesConcentration},
together with Lemma \ref{lema7} , the weak formulation associated to (\ref{AuxilaryProblemMalpha}) becomes:
find $\eta_\alpha \in Y$ such that
\[
\theta (  \nabla m_\alpha \, \nabla \phi ) - U \left(m_\alpha , \frac{\partial \phi}{\partial x_3}\right)
= 0 .
\]

By using Lax-Milgram lemma and standard computations when $U$ is small enough satisfying (\ref{SmallnessConditionOnUAuxiliaryProblem}),
recalling also (\ref{ChangeOFVariablesConcentrationAuxiliary}),
we obtain that there  is a unique solution $\eta_\alpha$ of the previous weak formulation and, moreover,
\begin{equation}
\label{EstimateNablaMalpha}
|\nabla m_\alpha | = |\nabla \eta_\alpha |   
\leq U \frac{2 \alpha}{\theta |\Omega|^{1/2}}.
\end{equation}

Now, we recall that in $H^1 (\Omega)$, the expressions $\| \phi \|_1$ and $|\int_\Omega \phi \, dx | + | \nabla \phi |$ give equivalent norms.
In particular, this leads to
\[
| m_\alpha | \leq \| m_\alpha \|_1 
\leq 
C \, [  \int_\Omega  m_\alpha  \, dx | + | \nabla  m_\alpha  | \, ]
\leq 
C \, [ \, \alpha + U \frac{2 \alpha}{\theta |\Omega|^{1/2}} \, ] .
\]

Estimate $ | \Delta m_\alpha |  $ is obtained directly from the equation and (\ref{EstimateNablaMalpha}).
In fact, using the equation and then  (\ref{EstimateNablaMalpha}), we obtain that:
\[
 | \Delta m_\alpha |  \leq \frac{U}{\theta}  | \nabla m_\alpha| \leq \frac{2 U^2  \alpha}{\theta^2 |\Omega|^{1/2}}.
\]

Next, we observe that the just obtained $m_\alpha$ is a distributional solution of (\ref{AuxilaryProblemMalpha}); 
from Theorem 10.1, p. 188-1899, of Ladyzhenskaya, Ural'tseva \cite{LU}, we have $m_\alpha \in H^2 (\Omega)$  
and $\| m_\alpha \|_ 2 \leq C [ \, |m_\alpha | + | \Delta m_\alpha | \, ] $.

By using this and the previous results, we then get the stated estimate for $\| m_\alpha \|_ 2 $.

Finally, since Lemma \ref{ineq} implies that 
\[
\vert \nabla m_\alpha \vert_4  \leq \vert   \nabla m_\alpha \vert^{1/4}\Vert  \nabla m_\alpha \Vert_1^{3/4} 
\leq \vert   \nabla m_\alpha \vert^{1/4}\Vert  m_\alpha \Vert_2^{3/4}, 
\]

\noindent
the just obtained estimates imply that
$|\nabla m_\alpha |_4$ decreases to zero as $ U$ approaches zero.
\hfill $\Box$

\begin{remark}
We stress that in the present section and in Section \ref{GlobalEvolutionaryStrongSolutions}, $m_\alpha$ is always the solution of problem (\ref{AuxilaryProblemMalpha}) given by Lemma \ref{AuxiliaryProblemLemma}.

\end{remark}

From now on in the rest of this work, we assume that $U$ is at least small enough such that (\ref{SmallnessConditionOnUAuxiliaryProblem}) holds.
More stringent smallness conditions will be explicited when necessary.

Next, similarly as in the previous section, to prove the existence of evolutionary strong solutions, it is convenient to introduce the following change of variables
\begin{equation}
\label{ChangeOFVariablesConcentrationOther}
\eta= m- m_\alpha ,
\end{equation}

\noindent
which also rewrites \eqref{TotalMassConservation} as
\[
\int_\Omega \eta (x, t) dx = 0, \; \mbox{for all} \; t \in [0,T] .
\]

Also, owing to our choice of $m_\alpha$, $\eta$ satisfies the same boundary condition as $m$, that is,
\[
\theta \displaystyle\frac{\partial \eta}{\partial \n} - U\eta n_3 = 0,
\]

\noindent
which is exactly the boundary condition (\ref{bc0EigenvectorsConcentration}) associated to the eigenfunctions of operator $A_1$.

Similarly as in the previous section, due to our change of variables \eqref{ChangeOFVariablesConcentrationOther}, to prove Theorem  \ref{localsolMainResults}, it is enough to prove the following result.

\begin{proposition} 
\label{localsol}
Let $\f \in L^2 (0,T; \X)$,  $\b_1 = 0$,
$\u_0 \in \J_0 = Dom(A^{1/2})$ and $\eta_0 \in Dom(A_1 ) \subset H^2 (\Omega) \cap Y$. 
Assume also that $\nu$
is a function of class  $C^1$ satisfying (\ref{nus}) and
$\displaystyle\sup_{\mathbb{R}} |\nu'(m)| < \nu_1' < \infty$. 

Then, for $U$ sufficiently
small, there exists $T^* = T^* (\Omega, \u_0, \eta_0, \nu) \leq T$,
such that problem   \eqref{WeakFormulationOtherProposition} and \eqref{InitialConditionOtherProposition}, has a solution $(\u, \eta)$
satisfying
\begin{eqnarray}
&& \u \in L^\infty (0, T^*, \J_0)   \cap L^2 (0, T^*; \J_0 \cap \H^2 (\Omega)), \quad \displaystyle\frac{\partial \u}{\partial t} \in L^2 (0, T^*; \L^2 (\Omega)), \label{r21}\\
&& \eta \in L^\infty (0, T^*; H^2 (\Omega) \cap Y), \quad \displaystyle\frac{\partial \eta}{\partial t} \in L^\infty (0, T^*; L^2 (\Omega)),\label{r22}
\end{eqnarray}
\begin{equation}
\begin{array}{l}
\u(t) \rightarrow \u_0 \ \mbox{strongly in} \ \J_0~ \quad{\rm as} \ t \rightarrow 0^+,   \\
 \eta(t) \rightarrow\eta_0 \  \mbox{strongly in}  \ \H^1 (\Omega) \, \mbox{and weakly in} \ \H^2 (\Omega) \; \mbox{as} \ t \rightarrow 0^+,
\end{array}
\end{equation}

\noindent
and
\begin{equation} 
\label{sistema2Other}
\left\{\begin{array}{l} 
\displaystyle
P \left(
\partial_t \u
- 2 \ {\rm div}\, (\nu(\eta + m_\alpha) D(\u))
+ \u\cdot \nabla \u 
+ \left(\eta + m_\alpha \right) \chi - \f 
\right) 
= 0 \quad {\rm in} \ L^2(0,T^*;\X),
\vspace{0.1cm}
\\
\displaystyle
\overline{P} 
\left(
\partial_t \eta 
- \theta \Delta \eta 
+ \u\cdot \nabla \eta 
+ u\cdot \nabla  m_\alpha  
+ U \displaystyle\frac{\partial \eta}{\partial x_3} 
\right)
=
0 
\quad {\rm in} \ L^\infty(0, T^*;L^2 (\Omega)).
\end{array}\right.
\end{equation}

\noindent
Here $P$ and $\overline{P}$ are respectively the projections defined in Lemmas \ref{Lema 1.8}  and \ref{Lema 1.9}.

\end{proposition}

\subsection{Preparation for the proof of Proposition  \ref{localsol}}

Before starting the proof of Proposition \ref{localsol}, let us briefly describe two key ideas to be used in the proof.

The weak formulation, written in terms of $(\u, \eta)$, is the following:
\begin{equation}
\label{WeakFormulationEquationsOthers}
 \left\{
\begin{array}{l}
\displaystyle
(\partial_t \u, \w)
  + 2(\nu(\eta + m_\alpha )D(\u), D(\w)) 
   + (\u\cdot\nabla \u, \w^j)
 + (\eta \chi, \w^j)
 \\
 \hspace{2cm}
 \displaystyle
 =   -  (m_\alpha \chi, \w)
+ (\f , \w),
\vspace{0.1cm}
  \\
  \displaystyle
( \partial_t \eta, \phi)
  + \theta (\nabla\eta, \nabla \phi)
  + (\u\cdot\nabla \eta,  \phi)
  + (\u\cdot\nabla m_\alpha, \phi)
- U\left(\eta, \frac{\partial  \phi }{\partial x_3}\right) 
=    0 ,   
 \vspace{0.1cm}
 \\
 \u =  \u_0 \ , \quad \ \eta = \eta_0 \quad \mbox{at time} \; t=0,
\end{array}
\right.
\end{equation}

\noindent
for all $\w \in \J_0$ and $\phi \in B$.

Then, observing (\ref{WeakFormulationEquationsOthers}), to build Galerkin approximations,  one would usually start by considering the weak formulation 
associated to (\ref{sistema2Other}); that is,
\begin{equation}
\label{ApproximateEquationsOthers}
 \left\{
\begin{array}{l}
\displaystyle
(\partial_t \u^n, \w^j)
  + 2(\nu(\eta^n + m_\alpha )D(\u^n), D(\w^j)) 
   + (\u^n\cdot\nabla \u^n, \w^j)
 + (\eta^n \chi, \w^j)
 \\
 \hspace{2cm}
 \displaystyle
 =   -  (m_\alpha \chi, \w^j )
+ (\f , \w^j),
\vspace{0.1cm}
  \\
  \displaystyle
( \partial_t \eta^n, \phi^\ell)
  + \theta (\nabla\eta^n, \nabla \phi^\ell)
  + (\u^n\cdot\nabla \eta^n,  \phi^\ell)
  + (\u^n\cdot\nabla m_\alpha, \phi^\ell)
- U\left(\eta^n, \frac{\partial  \phi^\ell }{\partial x_3}\right) 
=    0 ,   
 \vspace{0.1cm}
 \\
 \u^n(0) = \u_0^n , \quad  \eta^n(0) = \eta_0^n ,
\end{array}
\right.
\end{equation}

\noindent
for $1 \leq j, \ell \leq n$ and $\w^j$ and $\phi^\ell$ respectively in suitable Schauder bases $\{\w^i\}_{i=1}^\infty$ and $\{ \phi^k \}_{k=1}^\infty$,
with suitable initial conditions $\u_0$ and $\eta_0^n$ to be described later on.

However, since we will use an argument of smallness combined with an argument of contradiction, 
it is convenient to do a further change of variable related to the Galerkin approximations of the microorganism concentration in the equations  (\ref{ApproximateEquations}).
We define a new variable $\xi^n$ by:
\begin{equation}
\label{FurtherChangeOfVariableConcentration}
\xi^n = \eta^n - \eta^n_0.
\end{equation}

Since by this construction $\xi^n(0) = 0$, the time continuity implies that $\xi^n(t)$ will be small for small times $t$,
and this fact will be basic to reach a contradiction that will prove the proposition.

Another important point is that, due to strong nonlinearity of the second operator in the first equation, 
to obtain higher order estimates for $|D(\u^n)|$ and $|\nabla\xi^n|$, we will use the Helmholtz decomposition of $-\Delta\u$,
see Lemma \ref{helmholtz}. 
This introduces a auxiliary pressure-like term, whose estimation depends on  $|\nabla\xi^n|_4$. 
We will choose $T^*> 0$ sufficiently small  such that $|\nabla \xi^n|_{L^\infty(0.T^*; L^4 (\Omega))}$ is bounded independently of $n$.

We rewrite \eqref{ApproximateEquationsOthers} using the change of variable \eqref{FurtherChangeOfVariableConcentration} and
recalling that in the present case $\b_1 = 0$.
Moreover, instead of using  Galerkin approximations built with general Schauder bases, as we did in the proof of the previous proposition,
here and in the rest of this paper, we use spectral bases, which exist from the Lemmas \ref{StokesOperatorEigenvaluesEigenvectors}
and \ref{LaplacianEigenvaluesEigenvectors}.
More specifically, we use the eigenfunctions $\{ \w^i \}_{i=1}^\infty$ of the Stokes operator $A$ and the eigenfunctions $\{ \phi^k \}_{k=1}^\infty$ of operator $A_1$ as Schauder bases of  $\J_0$ and $B$, respectively.
Thus, we build approximate solutions by using the following Galerkin approximations:

\begin{equation}
\label{galerkinsol1}
\u^n(t,x) = \sum^n_{j=1} c_{n,j} (t)  \w^j(x);
\quad 
\xi^n (t,x) = \sum^n_{\ell=1} d_{n,\ell} (t)  \phi^\ell(x),
\end{equation}

Then, $\u^n$ and $\xi^n$ must satisfy the following approximate problem: for all $1 \leq j, \ell \leq n$,
\begin{equation}
\label{ApproximateEquationsStrong}
 \left\{
\begin{array}{l}
\displaystyle
(\partial_t \u^n,  \w^j)
  + 2(\nu(\xi^n  + \eta^n_0 +  m_\alpha )D(\u^n), D( \w^j))
   + (\u^n\cdot\nabla \u^n,  \w^j)
 + ( \xi^n   \chi,  \w^j)
 \\
 \hspace{2cm}
 \displaystyle
 =   
 - (  \eta^n_0  \chi,  \w^j)
 -  ( m_\alpha \chi, \w^j ) 
+ (\f , \w^j)
,
\vspace{0.1cm}
  \\
  \displaystyle
( \partial_t \xi^n, \phi^\ell)
  + \theta (\nabla\xi^n, \nabla  \phi^\ell)
  + \theta (\nabla\eta_0^n, \nabla  \phi^\ell)
  + (\u^n\cdot\nabla \xi^n, \phi^\ell)
  + (\u^n\cdot\nabla \eta^n_0, \phi^\ell) 
\vspace{0.2cm}
\\
\hspace{2cm}
\displaystyle
 + (\u^n\cdot\nabla m_\alpha, \phi^\ell) 
- U\left( \xi^n  , \frac{\partial  \phi^\ell}{\partial x_3}  \right) 
=
U\left( \eta^n_0  , \frac{\partial \phi^\ell}{\partial x_3} \right)  ,   
 \vspace{0.1cm}
 \\
 \u^n(0) = P_n \u_0 \ , \ \xi^n(0) = 0,
\end{array}
\right.
\end{equation}

Moreover, we take the initial conditions for the spectral approximations as follows:
\begin{equation}
\label{InitialConditionsProjections}
\u_0^n = P_n \u_0 \ , \quad \eta_0^n = \overline P_n \eta_0 ,
\end{equation}

\noindent
where $P_n$ and $\overline{P}_n$ are the projections described in the commentaries after respectively Lemma \ref{StokesOperatorEigenvaluesEigenvectors} and Lemma \ref{Lema 1.9}. 
Then, by using the properties described in those commentaries,
for the initial conditions as stated in Proposition \ref{localsol}, we have that
\begin{equation}
\label{InitialDataConvergences1}
\begin{array}{l}
|\u_0^n | \leq |\u_0 |, \quad |A^{1/2}\u_0^n | \leq |A^{1/2} \u_0 |,
\vspace{0.2cm}
\\
\u_0^n \rightarrow \u_0 \quad \mbox{in the norm of} \; \J_0 \; \mbox{and so, in the norm of} \;  \H^1(\Omega).
\end{array}
\end{equation}
\begin{equation}
\label{InitialDataConvergences2}
\begin{array}{l}
|\eta_0^n | \leq |\eta_0 |, \quad |A_1^{1/2}\eta_0^n | \leq |A_1^{1/2} \eta_0 |, \quad |A_1 \eta_0^n | \leq |A_1 \eta_0 |,
\vspace{0.2cm}
\\
\eta_0^n \rightarrow \eta_0 \quad \mbox{in the norm } \; |A (\cdot) | \; \mbox{and so, in the norm of} \; H^2(\Omega) .
\end{array}
\end{equation}

\subsection{Higher order estimates}

\begin{lemma}
For any given $\epsilon >0$, there is a positive constant $C_\epsilon$ such that
\begin{equation}
\label{332}
\begin{array}{c}
\displaystyle
\frac{d}{dt}|D(\u^n)|^2 
  + \frac{\nu_0}{2} | A\u^n |^2
  \vspace{0.2cm}
  \\
  \displaystyle
  \leq
 C |D(\u^n)|^6 
+ C_\eps  \left( |\nabla\xi^n|^8_4 
 + |\nabla m_\alpha|_4^8 +  |\nabla m_\alpha|_4^2 
+ C_{\eta_0}^8  +C_{\eta_0}^2 \right) |D(\u^n)|^2  
  \vspace{0.2cm}
  \\
  \displaystyle
+  C |\xi^n|^2
+ C C_{\eta_0}^2 
+C  |m_\alpha|^2 
+ C |\f|^2 
+ \eps \bar{C} \left( |\nabla\xi^n|_4 + C_{\eta_0} \right) |A\u^n|^2 .
\end{array}
\end{equation}

\end{lemma}

\noindent
{\bf Proof:}
By multiplying the first equation in \eqref{ApproximateEquationsStrong} by $\alpha_j c_{n,j}$,
 adding the result for $j=1, \ldots, n$,
 recalling that $A \u^n = \Sigma_{j=1}^n \alpha_j c_{n,j} \w^j$,
 and performing an integration by parts,
 we obtain the following identity
\begin{equation}
\label{StrongEstimateVelocity1}
\begin{array}{c}
\displaystyle
(\partial_t \u^n,  A\u^n)
  - 2 ( \, {\rm div} (\nu(\xi^n  + \eta^n_0 +  m_\alpha )D(\u^n) ), A\u^n \, )
   + (\u^n\cdot\nabla \u^n,  A\u^n)
 + ( \xi^n   \chi,  A\u^n)
 \\
 \displaystyle
 =   
 - (  \eta^n_0  \chi,  A\u^n)
 -   ( m_\alpha \chi, A\u^n)
+ (\f , A\u^n) .
\end{array}
\end{equation}

We can write the first term of the identity \eqref{StrongEstimateVelocity1} as 
\begin{equation}
\label{StrongEstimateVelocity1FistTerm}
\displaystyle
(\partial_t \u^n, A\u^n)=\frac{d}{dt}|D(\u^n)|^2 - 2\int_\Gamma \partial_t \u^n \cdot D(\u^n)\n~dS.
\end{equation}

\noindent
By splitting $\partial_t \u^n$ and $D(\u^n)\n$ in their tangential and normal components and taking into account the boundary conditions, we obtain that
\begin{equation}
\label{StrongEstimateVelocity1FistTermBoundary}
\begin{array}{ll}
\jnt_\Gamma \partial_t \u^n \cdot D(\u^n)\n~dS & = \jnt_\Gamma (\partial_t \u^n \cdot\n) ((D(\u^n)\n)\cdot\n)~dS
\\
&  +\jnt_\Gamma (\partial_t \u^n-(\partial_t \u^n\cdot\n)\n)\cdot [D(\u^n)\n -  \n \cdot(D(\u^n)\n) \n ]~dS
\\
& =0,
\end{array}
\end{equation}

\noindent
since $\partial_t \u^n \cdot\n =  \sum^n_{j=1} c_{n,j}' (t) (  \w^j \cdot\n )= 0 $ from the boundary condition (\ref{bc0EigenvectorsVelocity}),
and also
 $[D(\u^n)\n-  \n \cdot (D(\u^n)\n )\n ] = ( \sum^n_{j=1} c_{n,j} (t)\, [ D( \w^j )\n-\n \cdot( D( \w^j )\n )\n) ] = 0$ due again to boundary condition  (\ref{bc0EigenvectorsVelocity}) when $\b_1 = \mathbf{0}$.

Next, by using \eqref{StrongEstimateVelocity1FistTerm}, \eqref{StrongEstimateVelocity1FistTermBoundary}
and also the identity 
${\rm div}\, (\nu(m)D(\u)) =  \nu(m) \frac{1}{2}\Delta \u + \nu'(m) D(\u) \nabla m$,
which holds when ${\rm div} \, \u =0$, 
and the Helmholtz decomposition: $- \Delta \u^n = A\u^n + \nabla \overline q^n$ (see Lemma~\ref{helmholtz}),   \eqref{StrongEstimateVelocity1} can be rewritten as
\begin{equation}
\label{StrongEstimateVelocity2}
\begin{array}{c}
\displaystyle
\frac{d}{dt}|D(\u^n)|^2 
  + ( \,  \nu(\xi^n  + \eta^n_0 +  m_\alpha )  A\u^n   , A\u^n \,  )
  =
  - ( \,  \nu(\xi^n  + \eta^n_0 +  m_\alpha )  \nabla \overline q^n  , A\u^n \,  )
  \vspace{0.1cm}
  \\
  \displaystyle
    + 2 ( \,  \nu^\prime (\xi^n  + \eta^n_0 +  m_\alpha )  D (\u^n) (\nabla \xi^n + \nabla \eta^n_0 + \nabla m_\alpha)  , A\u^n \, )
    \vspace{0.1cm}
  \\
   - (\u^n\cdot\nabla \u^n,  A\u^n)
 -  (  (\xi^n  +  \eta^n_0  + m_\alpha)  \chi, A\u^n) 
+ (\f , A\u^n) 
  \vspace{0.1cm}
 \\
 \displaystyle
= \Sigma_{i=1}^5 I_i .
\end{array}
\end{equation}

By using the H\"older and Korn inequalities, Sobolev embedding,
Young inequality and (\ref{nus}), the terms in the right-hand side can be estimated
as follows. 

Since ${\rm div} A\u^n = {\rm div} A ( \sum^n_{j=1} c_{n,j} (t)  \w^j)  =  \sum^n_{j=1} c_{n,j} (t) \alpha^j {\rm div } \, \w^j= 0$,  we get
\begin{eqnarray*}
I_1 =
-(\nu(\xi^n+\eta^n_0+ m_\alpha)\nabla\overline q^n, A\u^n) & = &
(\overline q^n, \nabla (\nu(\xi^n + \eta^n_0 + m_\alpha) \cdot A\u^n)) \\
& = & (\overline q^n, \nu'(\xi^n+\eta^n_0 + m_\alpha) \nabla (\xi^n + \eta^n_0 + m_\alpha ) \cdot A\u^n).
\end{eqnarray*}
\noindent
By using the estimates of Lemma~\ref{helmholtz}, we get that
\begin{eqnarray*}
|I_1|&\leq& 
 \nu_1'|\overline q^n|_4 |\nabla (\xi^n+\eta^n_0 + m_\alpha)|_4 |A\u^n|
\\
&\leq & C|\overline q^n|^{1/4} \|\overline q^n\|_1^{3/4} |\nabla (\xi^n+\eta^n_0 + m_\alpha)|_4 |A\u^n| 
\\
&\leq & C|\nabla (\xi^n+\eta^n_0 + m_\alpha)|_4  \left[ \, C_\epsilon |\nabla \u^n|  + \epsilon |A \u^n| \, \right]^{1/4} |A \u^n |^{3/4} |A\u^n| 
\\
&\leq & C|\nabla (\xi^n+\eta^n_0 + m_\alpha)|_4  C \left[ \, C_\epsilon^{1/4} |\nabla \u^n|^{1/4}  + \epsilon^{1/4} |A \u^n|^{1/4} \, \right] |A \u^n |^{7/4} 
\\
&\leq&
C_\eps|\nabla (\xi^n+\eta^n_0 + m_\alpha)|_4 |\nabla\u^n|^{1/4}|A\u^n|^{7/4}+
\eps |\nabla (\xi^n+\eta^n_0 + m_\alpha )|_4  |A\u^n|^2
\\
&\leq& C_{\eps, \delta} |\nabla(\xi^n + \eta^n_0+ m_\alpha) |^8_4 |D(\u^n)|^2 + \delta |A\u^n|^2 + \eps  |\nabla(\xi^n + \eta^n_0 + m_\alpha)|_4 |A\u^n|^2 ,
\end{eqnarray*}

\noindent
with any $\varepsilon = \epsilon^{1/4}> 0$.

\begin{eqnarray*}
|I_2|&=& |2(\nu'(\xi^n+\eta^n_0+  m_\alpha  ) D(\u^n) (\nabla ( \xi^n + \eta^n_0 +  m_\alpha) ), A\u^n)| 
\\
&& \leq 2\nu_1' |D(\u^n)|_4   |\nabla( \xi^n +  \eta^n_0 + m_\alpha ) |_4  |A\u^n| 
\\
&& \leq 2\nu_1' |\nabla ( \xi^n +  \eta^n_0 + m_\alpha) |_4  |D(\u^n)|^{1/4} |A\u^n|^{7/4} 
\\
&& \leq C_\delta  |\nabla ( \xi^n +  \eta^n_0 + m_\alpha ) |_4^2 |D(\u^n)|^2  + \delta |A\u^n|^2 .
\end{eqnarray*}

\[
\begin{array}{l}
|I_3|=|(\u^n\cdot\nabla \u^n, A\u^n)|  
\leq  |\u^n|_6 |\nabla \u^n|_3 |A\u^n| 
\vspace{0.1cm}
\\
\hspace{0.7cm}
 \leq  |\nabla \u^n| \ |\nabla \u^n|^{1/2} |A\u^n|^{3/2} 
 \leq  C_\delta |D(\u^n)|^6 + \delta |A\u^n|^2.
\end{array}
\]

\begin{eqnarray*}
|I_4|=| ( (\xi^n+\eta^n_0 + m_\alpha)  \chi, A\u^n) | & \leq & C_\delta |\xi^n+\eta^n_0  + m_\alpha |^2 + \delta |A\u^n|^2.
\end{eqnarray*}

\[
|I_5| = | (\f , A\u^n)|=C_\delta |\f|^2 + \delta |A\u^n|^2 .
\]

By using these previous estimates in (\ref{StrongEstimateVelocity2}) and fixing $\delta$ sufficiently small, we obtain that
\begin{equation}
\label{332Preparatory}
\begin{array}{c}
\displaystyle
\frac{d}{dt}|D(\u^n)|^2 
  + \frac{\nu_0}{2} | A\u^n |^2
  \vspace{0.2cm}
  \\
  \displaystyle
  \leq
  \displaystyle
C_\eps |\nabla(\xi^n + \eta^n_0  + m_\alpha)|^8_4 |D(\u^n)|^2  
    + C |\nabla (\xi^n+\eta^n_0 +  m_\alpha)|^2_4 |D(\u^n)|^2  
+ C |D(\u^n)|^6 
  \vspace{0.2cm}
  \\
  \displaystyle
+  C |\xi^n+\eta^n_0  +  m_\alpha |^2 
+ C |\f|^2 
+ \eps  |\nabla(\xi^n + \eta^n_0 +  m_\alpha )|_4 |A\u^n|^2 .
\end{array}
\end{equation}

Now, by using estimates (\ref{EstimateImportant}), we observe that
\begin{equation}
\label{EstimateEta0n}
\max \{  |\nabla \eta^n_0|,
 |\nabla \eta^n_0|_4 \}
\leq 
C |A_1 \eta_0^n |
\leq 
C |A_1 \eta_0|
= C_{\eta_0},
\end{equation}

\noindent
which, together with (\ref{332Preparatory}), implies (\ref{332}).
\hfill $\Box$

\begin{lemma}
\label{EstimatePartialunt}
\begin{equation}
\label{336}
\begin{array}{l}
\displaystyle
|\partial_t \u^n |^2
\leq   C | A \u^n |^2 +   C |  \nabla   \xi^n|_4^2 
+   C C_{\eta_0}
 + 24 |m_\alpha |^2 
+24 | \f |^2
\end{array}
\end{equation}

\end{lemma}

\noindent
{\bf Proof:}
By taking
 $\partial_t\u^n (\cdot, t)$ as test function in the first equation of   (\ref{ApproximateEquationsStrong}),  
using Lemma \ref{lema7} and the fact that $\b_1 = 0$,
we obtain that
\begin{equation}
\label{vt}
\begin{array}{l}
\displaystyle
|\partial_t \u^n |^2
=
  2 (  {\rm div} \left(   \nu(\xi^n  + \eta^n_0 +  m_\alpha )D(\u^n)  \right), \partial_t \u^n )
   - (\u^n\cdot\nabla \u^n,  \partial_t \u^n)
 - ( \xi^n   \chi,  \partial_t \u^n)
 \\
 \hspace{2cm}
 \displaystyle
 - (  \eta^n_0  \chi,  \partial_t \u^n)
 -   (   m_\alpha  \chi, \partial_t \u^n)
+ (\f , \partial_t \u^n) ,
\end{array}
\end{equation}

\noindent
which implies
\[
\begin{array}{l}
\displaystyle
|\partial_t \u^n |^2
\leq
   24 |  {\rm div} \left(   \nu(\xi^n  + \eta^n_0 +    m_\alpha )D(\u^n)  \right) |^2
   + 24  | \u^n\cdot\nabla \u^n |^2
   + 24  | \xi^n |^2
 \\
 \hspace{2cm}
 \displaystyle  
 + 24 | \eta^n_0 |^2
 + 24 |  m_\alpha|^2
+ 24 | \f |^2 .
\end{array}
\]

By using  Young's inequality and Lemma \ref{ineq}, from this last result we get that
\[
\begin{array}{l}
\displaystyle
|\partial_t \u^n |^2
\leq
  C | A \u^n |^2
+ C  | D(\u^n)  |_4^2
   + C | \u^n|_4^2 + C |\nabla \u^n |_4^2
+   C |  \nabla   \xi^n|_4^2 
 \\
 \hspace{1.6cm}
 \displaystyle  
   + 24  | \xi^n |^2
 + 24 | \eta^n_0 |^2
+   C |  \nabla  \eta^n_0|_4^2
 + 24 |m_\alpha |^2
+24 | \f |^2
 \\
 \hspace{1.4cm}
 \displaystyle  
\leq   C | A \u^n |^2 +   C |  \nabla   \xi^n|_4^2 
   + 24  | \xi^n |^2
+   C |  \nabla  \eta^n_0|_4^2
 + 24  |m_\alpha |^2
+24 | \f |^2
 \\
 \hspace{1.4cm}
 \displaystyle  
\leq   C | A \u^n |^2 +   C |  \nabla   \xi^n|_4^2 
   + 24  | \xi^n |^2
+   C C_{\eta_0}
 + 24  |m_\alpha |^2
+24 | \f |^2 .
\end{array}
\]

\noindent
where we used estimate (\ref{EstimateEta0n}) to obtain the last line.

By recalling that  $\xi^n$ has mean value zero, we have   $ | \xi^n | \leq C |\nabla \xi^n | \leq C |\nabla \xi^n |_4 $, with constants independent of $n$,
and the last estimate implies (\ref{336}).
\hfill $\Box$

\begin{lemma}
\begin{equation}
\label{334}
 |\nabla\xi^n (t) |^2 + \int^t_0 |A_1 \xi^n (s)|^2 ds   
\leq 
H_1 (t;  \|D(\u^n)\|_{L^\infty(0, t; \L^2 (\Omega))}),
\end{equation}

\noindent
where $H_1 (t; z) \equiv 2 C C_{\eta_0}^2 t \exp \left(M_1(z) \, t \right) $  and 
$M_1 (z) = C (z^4 + 1)$.

\noindent
We stress that function $ H_1 $ is strictly increasing in each one of its variables and is independent of $n$.

\end{lemma}

\noindent
{\bf Proof:}
By multiplying the second equation in \eqref{ApproximateEquationsStrong} by $- \beta_\ell d_{n, \ell}$,
 adding the result for $\ell=1, \ldots, n$,
 recalling that $-A_1 \xi^n = - \Sigma_{\ell=1}^n \beta_j d_{n, \ell} \phi^\ell$,
 and performing an integration by parts,
 we obtain the following:
 \begin{equation}
 \label{StrongEstimateConcentration1}
 \begin{array}{l}
 \displaystyle
- ( \partial_t \xi^n, A_1 \xi^n)
  + \theta | A_1\xi^n |^2
  =     - \theta ( A_1\eta_0^n, A_1 \xi^n)
   - \theta ( \Delta m_\alpha, A_1 \xi^n)
\vspace{0.1cm}
  \\
  \displaystyle
  +   (\u^n\cdot\nabla \xi^n, A_1 \xi^n)
  +   (\u^n\cdot\nabla  \eta_0^n, A_1 \xi^n)
  +   (\u^n\cdot\nabla m_\alpha, A_1 \xi^n)
\vspace{0.1cm}
  \\
  \displaystyle
+ U \left(\frac{\partial}{\partial x_3} (\xi^n + \eta^n_0 ) , A_1 \xi^n\right).
\end{array}
 \end{equation}

Integrating by parts the first term of \eqref{StrongEstimateConcentration1}  and taking into account that the fourth boundary condition in (\ref{bc0}) is true for $\eta^n_0 = \overline P_n \eta_0$, we have
\begin{eqnarray}
-(\partial_t\xi^n, A_1\xi^n) & = & (\nabla\partial_t \xi^n, \nabla \xi^n) - \int_{\partial\Omega} \partial_t \xi^n \frac{\partial \xi^n}{\partial n} dS
\nonumber
\\
& = & \frac12 \frac{d}{dt} |\nabla \xi^n|^2 - \frac{U}\theta \int_{\partial\Omega} \partial_t \xi^n\xi^n n_3 dS 
\nonumber
\\
& = & \frac12 \frac{d}{dt} |\nabla \xi^n|^2 - \frac12 \frac{U}\theta \frac{d}{dt} \int_{\partial \Omega} |\xi^n|^2 n_3 dS. 
\nonumber
\\
& = &
 \frac12 \frac{d}{dt} \left(|\nabla \xi^n|^2 - \frac{U}{\theta}
\int_{\partial\Omega} |\xi^n|^2 n_3 dS\right).
\label{ReescritaDaDerivadaTemporal}
\end{eqnarray}

We also remark that, denoting
\begin{equation}
\label{yn(t)}
y^2_n(t) 
= 
|\nabla\xi^n|^2 - \displaystyle\frac{U}{\theta} \int_{\partial\Omega} |\xi^n|^2 n_3 dS ,
\end{equation}

\noindent
we immediately have $y^2_n(t) \leq |\nabla \xi^n|^2 $. 
On the other hand, we also have that
$ \displaystyle y^2_n(t)  \geq |\nabla\xi^n|^2 - \displaystyle\frac{U}{\theta} C | \nabla \xi^n |^2 \geq \frac{1}{2}  | \nabla \xi^n |^2 $
if $U$ is small enough such that
\[
1 - \frac{U}{\theta} C \geq \frac{1}{2},
\]

\noindent
where $C$ is the positive constant associated to the continuity of trace operator from $H^1 (\Omega) \cap Y$ to $L^2 (\partial \Omega)$.
For such $U$, we have that
\begin{equation}
\label{yn(t)EquivalentNomGradxi}
 \frac{1}{2} |\nabla\xi^n|^2 \leq y^2_n(t) \leq |\nabla \xi^n|^2 .
\end{equation}

\noindent
That is, $|\nabla\xi^n|$ and $y_n$ are equivalent norms.

The terms on the right side of \eqref{StrongEstimateConcentration1} can be estimated as follows.
\[
| \theta ( A_1\eta_0^n, A_1 \xi^n) | \leq  C_\delta |A_1\eta_0^n |^2 + \delta |A_1 \xi^n |^2 ,
\]
\[
| \theta ( \Delta m_\alpha, A_1 \xi^n) | \leq  C_\delta \| m_\alpha\|_2^2 + \delta |A_1 \xi^n |^2 ,
\]
\[
\begin{array}{l}
 |(\u^n\cdot\nabla  \xi^n , A_1 \xi^n)| 
\leq 
|\u^n|_6 |\nabla \xi^n |_3 |A_1 \xi^n| 
\\
\hspace{3.1cm}
\leq 
C \Vert \u^n \Vert_1
\vert \nabla \xi^n  \vert^{1/2} \Vert \nabla \xi^n \Vert_1^{1/2}
|A_1 \xi^n| 
\vspace{0.2cm}
\\
\hspace{3.1cm}
 \leq 
C | D( \u^n ) | \,
| \nabla \xi^n  |^{1/2} |A_1 \xi^n|^{3/2}
\vspace{0.2cm}
\\
\hspace{3.1cm}
 \leq 
C_{\delta}( |D(\u^n)|^4
|\nabla \xi^n|^2
+ \delta |A_1 \xi^n|^2 ,
\end{array}
\]
\[
\begin{array}{l}
 |(\u^n\cdot\nabla   \eta_0^n , A_1 \xi^n)| 
\leq 
|\u^n|_4 |\nabla  \eta_0^n |_4 |A_1 \xi^n| 
\\
\hspace{3.1cm}
\leq 
C \Vert \u^n \Vert_1 |\nabla  \eta_0^n |_4
|A_1 \xi^n| 
\vspace{0.2cm}
\\
\hspace{3.1cm}
 \leq 
C_{\delta}( |D(\u^n)|^2
|\nabla  \eta_0^n |_4^2
+ \delta |A_1 \xi^n|^2 ,
\end{array}
\]
\[
\begin{array}{l}
 |(\u^n\cdot\nabla   m_\alpha , A_1 \xi^n)| 
\leq 
|\u^n|_6 |\nabla   m_\alpha  |_3 |A_1 \xi^n| 
\\
\hspace{3.1cm}
\leq 
C \Vert \u^n \Vert_1 \|   m_\alpha   \|_2
|A_1 \xi^n| 
\vspace{0.2cm}
\\
\hspace{3.1cm}
 \leq 
C_{\delta}( |D(\u^n)|^2
\|   m_\alpha  \|_2^2
+ \delta |A_1 \xi^n|^2 ,
\end{array}
\]
\[
|U \left(\frac{\partial}{\partial x_3} (\xi^n + \eta^n_0 ), A_1 \xi^n\right)| 
\leq 
C_\delta |\nabla\xi^n|^2 + C_\delta |\nabla \eta^n_0|^2   + \delta |A_1 \xi^n|^2.
\]

By using (\ref{ReescritaDaDerivadaTemporal}), (\ref{yn(t)}) and
the last estimates in (\ref{StrongEstimateConcentration1}), with a fixed appropriate small $\delta$,  it follows that

 \begin{equation}
 \label{paraxi}
 \begin{array}{l}
 \displaystyle
 \frac{1}{2} \frac{d}{dt} ( y_n (t)^2 )
  +  \frac{\theta}{2} | A_1\xi^n |^2
  \leq 
   C_{\delta}( |D(\u^n)|^4 |\nabla \xi^n|^2 
+ C_\delta |A_1\eta_0^n |^2 
 + C_\delta \| m_\alpha\|_2^2  
\vspace{0.2cm}
  \\
\hspace{4.3cm}
  \displaystyle 
+  C_{\delta}( |D(\u^n)|^2 \left(  |\nabla  \eta_0^n |_4^2  +  \|  m_\alpha  \|_2^2  \right) 
\vspace{0.2cm}
  \\
\hspace{4.3cm}
  \displaystyle
+ C_\delta |\nabla\xi^n|^2 + C_\delta |\nabla \eta^n_0|^2  .
\end{array}
 \end{equation}

By integrating this last inequality in time from $0$ to $t$, using (\ref{yn(t)EquivalentNomGradxi}) and (\ref{yn(t)}) at $t=0$,
recalling that $\xi^n (0) = 0$,
we obtain that
\[
\begin{array}{c}
\displaystyle
 C_1 |\nabla\xi^n|^2 + \int^t_0 |A_1 \xi^n|^2 ds  
\leq 
|y_n (t)|^2 + \int^t_0 |A_1 \xi^n|^2 ds  
\vspace{0.2cm}
\\
\displaystyle
\leq 
 C_\delta \left( |A_1\eta_0^n |^2 + |\nabla \eta^n_0|^2  + \| m_\alpha\|_2^2   \right) t
+ \int_0^t C_{\delta}( |D(\u^n)|^2 \left(  |\nabla  \eta_0^n |_4^2  +  \|  m_\alpha  \|_2^2  \right) ds 
\vspace{0.2cm}
\\
\displaystyle
+ \int^t_0 C_{\delta}\left( |D(\u^n)|^4 + 1 \right) |\nabla \xi^n|^2 ds .
\end{array}
\]

By recalling (\ref{EstimateEta0n}) and (\ref{MalphaEstimates}), the last inequality implies that
\[
\displaystyle
 |\nabla\xi^n|^2 + \int^t_0 |A_1 \xi^n|^2 ds  
\leq 
2 C C_{\eta_0 m_\alpha}^2 t
+ \int_0^t C |D(\u^n)|^2ds 
+ \int^t_0 C \left( |D(\u^n)|^4 + 1 \right) |\nabla \xi^n|^2 ds ,
\]

\noindent
which, by using Gronwall lemma, implies that
\[
\begin{array}{c}
\displaystyle
 |\nabla\xi^n (t) |^2 + \int^t_0 |A_1 \xi^n (s)|^2 ds  
\vspace{0.2cm}
\\
\displaystyle
\leq 
\left(2 C C_{\eta_0m_\alpha}^2 t    
+ \int_0^t C |D(\u^n)|^2ds 
\right)
\exp \left( \int^t_0 C \left( |D(\u^n)|^4 + 1 \right)  dt \right),
\end{array}
\]

\noindent
and so,
\[
\begin{array}{c}
\displaystyle
 |\nabla\xi^n (t) |^2 + \int^t_0 |A_1 \xi^n (s)|^2 ds  
\leq 
\vspace{0.2cm}
\\
\displaystyle
C t
(2  C_{\eta_0m_\alpha}^2    +  \|D(\u^n)\|_{L^\infty(0, t; \L^2 (\Omega))}^2  )
\exp \left(  C t ( \|D(\u^n)\|_{L^\infty(0, t; \L^2 (\Omega))}^4 + 1 ) \right),
\end{array}
\]

\noindent
from which we obtain (\ref{334}).
\hfill $\Box$

\begin{lemma}
There exist non-negative functions  $H_2$, $H_3$,
which are independent of $ n $ and strictly increasing in each one
of their variables, such that:
\begin{equation}
\label{parttxi}
|\partial_t \xi^n(t)|^2 
\leq  
H_2(t; \|\nabla\xi^n\|_{L^\infty(0,t;L^4)}; \|A\u^n\|_{L^2(0,t;\L^2)}).
\end{equation}
\noindent
and
\begin{equation}
\label{339}
|A_1\xi^n|^2 
\leq 
H_3 (t; \|D(\u^n)\|_{L^\infty(0, t; L^2)} ;
\|\nabla\xi^n\|_{L^\infty(0, t; L^4)} ; \|A\u^n\|_{L^2(0, t; \L^2)}).
\end{equation}

\end{lemma}

\noindent
{\bf Proof:}
To obtain estimates for $\partial_t \xi^n$, we differentiate the second equation of  (\ref{ApproximateEquationsStrong}) with respect to $t$ and take
$\partial_t \xi^n (\cdot, t)$ as test function to obtain that
\[
\begin{array}{l}
\displaystyle
 \frac12 \frac{d}{dt} |\partial_t\xi^n|^2 
+ \theta |\nabla\partial_t\xi^n|^2 
= 
- (\partial_t \u^n\cdot\nabla \xi^n , \partial_t \xi^n)
  - (\partial_t \u^n\cdot\nabla \eta^n_0, \partial_t \xi^n) 
\vspace{0.2cm}
\\
\displaystyle
\hspace{4.4cm}
 - (\partial_t\u^n\cdot\nabla m_\alpha, \partial_t \xi^n) 
+ U(\partial_t \xi^n, \frac{\partial}{\partial x_3} \partial_t \xi^n) 
\vspace{0.2cm}
\\
\displaystyle
\hspace{4.cm}
  \leq 
|\partial_t \u^n| \ |\nabla \xi^n |_4 |\partial_t \xi^n|_4 
 + |\partial_t \u^n| \ |\nabla \eta_0^n |_4 |\partial_t \xi^n|_4 

\vspace{0.2cm}
\\
\hspace{4.4cm}
\displaystyle
+|\partial_t \u^n| \ |\nabla m_\alpha |_4 |\partial_t \xi^n|_4 
+ U |\partial_t\xi^n|
\ |\nabla\partial_t \xi^n| 
\vspace{0.2cm}
\\
\hspace{4.cm}
\displaystyle
  \leq 
|\partial_t \u^n| \ |\nabla \xi^n |_4 \| \partial_t \xi^n \|_1
 + |\partial_t \u^n|  |\nabla \eta_0^n |_4 \| \partial_t \xi^n \|_1
\vspace{0.2cm}
\\
\hspace{4.4cm}
\displaystyle
 + |\partial_t \u^n \ |\nabla m_\alpha |_4 \| \partial_t \xi^n| \|_1
+ U |\partial_t\xi^n|
\ |\nabla\partial_t \xi^n|  
\vspace{0.2cm}
\\
\hspace{4.cm}
\displaystyle
  \leq 
C |\partial_t \u^n \|\nabla \xi^n |_4 | \nabla \partial_t \xi^n| 
 + C |\partial_t \u^n| |\nabla \eta_0^n |_4 | \nabla \partial_t \xi^n| 
\vspace{0.2cm}
\\
\displaystyle
\hspace{4.4cm}
 + C |\partial_t \u^n| |\nabla m_\alpha |_4 | \nabla \partial_t \xi^n| 
+ U |\partial_t\xi^n|
\ |\nabla\partial_t \xi^n|  
\vspace{0.2cm}
\\
\hspace{4.cm}
\displaystyle
 \leq C_\delta |\partial_t \u^n|^2 |\nabla \xi^n |^2_4 \,  
 + C_\delta |\partial_t \u^n|^2 |\nabla \eta_0^n |^2_4 
\vspace{0.2cm}
\\
\displaystyle
\hspace{4.4cm}
 + C_\delta |\partial_t \u^n|^2 |\nabla m_\alpha |^2_4 \,
+ \delta  |\nabla \partial_t \xi^n|_4^2
\vspace{0.2cm}
\\
\displaystyle
\hspace{4.4cm}
+ C_\delta |\partial_t \xi^n|^2 + \delta |\nabla\partial_t
\xi^n|^2 .
\end{array}
\]

By fixing an appropriate $\delta$ and recalling (\ref{EstimateEta0n}) and (\ref{MalphaEstimates}), we obtain:
\begin{equation}
\label{335}
\frac{d}{dt} |\partial_t \xi^n|^2 
+ \theta |\nabla\partial_t \xi^n|^2 
\leq 
C_{\eta_0m_\alpha} |\partial_t \u^n|^2  \left( 1 + |\nabla\xi^n|^2_4 \right)
+C  |\partial_t \xi^n|^2.
\end{equation}

Next, by taking
 $\partial_t \xi^n (\cdot, t)$ as test function in the second equation of   (\ref{ApproximateEquationsStrong}),   using Lemma \ref{lema7},
and proceeding as done in Lemma \ref{EstimatePartialunt},
we obtain that
\[
\begin{array}{l}
  \displaystyle
| \partial_t \xi^n |^2
=
   \theta (\Delta \xi^n,    \partial_t \xi^n)
  + \theta (\Delta \eta_0^n,  \partial_t \xi^n)
\vspace{0.2cm}
\\
\hspace{1.5cm}
  - (\u^n\cdot\nabla \xi^n,  \partial_t \xi^n)
 - (\u^n\cdot\nabla \eta^n_0,  \partial_t \xi^n  ) 
 - (\u^n\cdot\nabla m_\alpha,  \partial_t \xi^n  ) 
\vspace{0.2cm}
\\
\hspace{1.5cm}
\displaystyle
- U\left(  \frac{\partial   \xi^n }{\partial x_3}  , \partial_t \xi^n \right) 
- U\left( \frac{\partial  \eta^n_0 }{\partial x_3}  ,  \partial_t \xi^n \right)    
  . 
\end{array}
\]

\noindent
which, by working in a standard way, using (\ref{EstimateEta0n}) and (\ref{MalphaEstimates}),
 implies that
\begin{equation}
\label{EstimateXit}
\begin{array}{l}
  \displaystyle
| \partial_t \xi^n |^2
\leq
   C |A_1 \xi^n|^2
  + C |\Delta \eta_0^n|^2

\vspace{0.2cm}
\\
\hspace{1.4cm}
 +  C |\u^n|_4^4 
+ C |\nabla \xi^n |_4^4
 \; + C |\nabla \eta_0^n |_4^4 
 \; + C |\nabla m_\alpha |_4^4 
\vspace{0.2cm}
\\
\hspace{1.4cm}
+ C | \nabla \xi^n |^2
+ C | \nabla \eta^n_0|^2
\vspace{0.2cm}
\\
\hspace{1.25cm}
\leq
   C |A_1 \xi^n|^2
 +  C |D(\u^n)|^4 
+ C |\nabla \xi^n |_4^4
+ C | \nabla \xi^n |^2
  + C C_{\eta_0m_\alpha}^2.
\end{array}
\end{equation}

By recalling that $\xi^n (0) = 0$, from (\ref{EstimateXit}) we have that
\begin{equation}
\label{338}
\begin{array}{l}
  \displaystyle
| \partial_t \xi^n (0)|^2
\leq
  C |D(\u^n (0) |^4 
  + C C_{\eta_0m_\alpha}^2,
\leq
  C |A^{1/2}\u_0| ^4 
  + C C_{\eta_0m_\alpha}^2,
\end{array}
\end{equation}

\noindent
since $|D(\u^n (0)) | \leq C |A^{1/2}\u^n (0)|\leq C |A^{1/2}\u_0| $.

By using the Gronwall lemma in (\ref{335}), we get that
\begin{equation}
\label{335A}
\begin{array}{c}
\displaystyle
 |\partial_t \xi^n (t)|^2 
+ \theta \int_0^t  |\nabla\partial_t \xi^n|^2 ds
\vspace{0.2cm}
\\
\displaystyle
\leq 
  \left(  C |A^{1/2}\u_0| ^4    + C C_{\eta_0m_\alpha}^2 
+C_{\eta_0m_\alpha} \left(  1 +  \|\nabla\xi^n\|_{L^\infty(0,t;L^4)}^2  \right)  \|\partial_t \u^n \|_{L^2(0,t;\L^2)}^2 \right) \exp (Ct) .
\end{array}
\end{equation}

Since  (\ref{336}) implies that
\begin{equation}
\label{336A}
\begin{array}{l}
\displaystyle
\|\partial_t \u^n \|_{L^2(0,t;\L^2)}^2
\leq   C \| A \u^n \|_{L^2(0,t;\L^2)}^2 
+   C t \|  \nabla   \xi^n\|_{L^\infty(0,t;L^4)}^2 
\vspace{0.2cm}
\\
\displaystyle
\hspace{3cm}
+   C t^{1/2}
 + C \| m_\alpha \|_{L^2(0,t;L^2)}
+C \| \f \|_{L^2(0,t;\L^2)}^2 ,
\end{array}
\end{equation}

\noindent
by plugging (\ref{336A}) in (\ref{335A}), we obtain (\ref{parttxi}).

\vspace{0.1cm}
Next, from (\ref{StrongEstimateConcentration1}) and Young inequality, 
by proceeding similarly as before, 
we obtain that
 \begin{equation}
 \begin{array}{c}
 \displaystyle
  | A_1\xi^n (t) |^2
  \leq
C \left(
 | \partial_t \xi^n (t)|^2
+ |\u^n|_4^4 
+ |\nabla \xi^n (t)|_4^4
 + |\nabla \eta_0^n |_4^4
 + |\nabla m_\alpha|_4^4
\right.
\vspace{0.2cm}
\\
\displaystyle
\left.
 + | A_1\eta_0^n|^2
 + |\Delta m_\alpha |^2 
+ |\nabla \xi^n|^2 
+ |\nabla \eta^n_0 |^2
\right)
\vspace{0.2cm}
\\
\displaystyle
  \leq
C \left(
 | \partial_t \xi^n (t) |^2
+ |D ( \u^n) |_{L^\infty(0, t; \L^2)}^4 
+ |\nabla \xi^n|_{L^\infty(0, t; L^4)}^4
+ |\nabla \xi^n|_{L^\infty(0, t; L^4 (\Omega))}^2 
\right.
\vspace{0.2cm}
\\
\left.
\displaystyle
 + | A_1\eta_0^n|^2
+ |\nabla \eta^n_0 |^2
 + |\nabla \eta_0^n |_4^4
 + | \Delta m_\alpha|^2 
 + |\nabla m_\alpha |_4^4 
\right) .
\end{array}
 \end{equation}

By plugging  (\ref{parttxi}) in this last inequality, we obtain (\ref{339}).
\hfill $\Box$

\begin{lemma}
There is $T^*$, with $0 < T^* \leq T$ and independent of $n$,
and nonnegative continuous functions $\gamma(t)$, $H_1 (t, z_1)$, $H_2(t, z_1,z_2)$, $H_3 (t, z_1, z_2, z_3) $, $H_5(t)$,
which are nondecreasing in each of their independent variables and also independent of $n$,
 such that, for $ t \in [0, T^*]$, there hold
\begin{equation}
\label{xi1}
\begin{array}{l}
\displaystyle
|D(\u^n(t)|^2 + \frac{\nu_0}4 \int^t_0 |A\u^n|^2 ds 
\leq 
\gamma (t)
\leq 
\gamma (T^*), 
\vspace{0.2cm}
\\
\displaystyle
|\nabla \xi^n (t)|^2 \leq H_1 (t, \gamma (t)) \leq H_1 (T^*, \gamma (T^*)) ,
\vspace{0.2cm}
\\
|A_1 \xi^n(t)|^2 \leq H_3 (t, \gamma(t), 1, 4 \gamma(t)/\nu_0)) \leq H_3 (T^*, \gamma(T^*), 1, 4 \gamma(T^*)/\nu_0)),
\vspace{0.2cm}
\\
|\nabla\xi^n(t)|_4 < 1,
\vspace{0.2cm}
\\
|\partial_t \xi^n(t)|^2 
\leq  
H_2(t; 1; 4 \gamma (t)/\nu_0) \leq H_2(T^*; 1; 4 \gamma (T^*)/\nu_0) ,
\vspace{0.2cm}
\\
\displaystyle
\int_0^t |\partial_t \u^n |^2 ds
\leq   H_5(t) \leq H_5 (T^*).
\end{array}
\end{equation}

\end{lemma}

\noindent
{\bf Proof:}
From (\ref{G}), we have that $|\eta^n(t)|^2 \leq G(t)/(4 C_{\nu_0}) $, which, taking into account the definition (\ref{FurtherChangeOfVariableConcentration}) implies that
$|\xi^n |^2= 2 |\eta^n| + 2 |\eta^n_0|^2 \leq   G(t)/(2 C_{\nu_0}) + 2 |\eta^n_0|^2$.
By using this in the differential inequality (\ref{332}), we obtain that
\begin{equation}
\label{332Z}
\begin{array}{c}
\displaystyle
\frac{d}{dt}|D(\u^n)|^2 
  + \frac{\nu_0}{2} | A\u^n |^2
  \vspace{0.2cm}
  \\
  \displaystyle
  \leq
C |D(\u^n)|^6 
+ C_\eps  \left( |\nabla\xi^n|^8_4 
 + |\nabla m_\alpha|_4^8 +  |\nabla m_\alpha|_4^2 
+ C_{\eta_0}^8  +C_{\eta_0}^2 \right) |D(\u^n)|^2  
  \vspace{0.2cm}
  \\
  \displaystyle
+  C  G(t)/(2 C_{\nu_0}) + 2 C |\eta^n_0|^2
+ C C_{\eta_0}^2 
+C  |m_\alpha|^2 
+ C |\f|^2 
+ \eps \bar{C} \left( |\nabla\xi^n|_4 + C_{\eta_0} \right) |A\u^n|^2 .
\end{array}
\end{equation}

\vspace{0.1cm}
Next, we fix  $\eps > 0$ such that 
\begin{equation}
\label{epsilon}
\eps \bar{C} (1 + C_{\eta_0}  ) \leq \frac{\nu_0}{4} ,
\end{equation}

\noindent
where $\bar{C}$ is the constant appearing in the last term of (\ref{332Z}).
This also fix the positive constant $C_\epsilon$ appearing  in the second term of the right-hand side of (\ref{332Z}).

Once such such $\epsilon$ is chosen,  we will prove with the help of Lemma \ref{lema110} that there is $T^*$  independent of $n$ such that  $|\nabla\xi^n (t)|_4<1$ for $t \in [0,T^*)$.

For this, we observe the expressions of the terms in (\ref{332Z}) and,  using the notation of  Lemma \ref{lema110}, define
\begin{equation}
\label{Lemma10Correspondence}
\begin{array}{l}
\varphi (t) = |D(\u^n) (t)|^2, \quad
\varphi_0 = |D(\u^n_0)|^2 
\vspace{0.2cm}
\\
\displaystyle
\psi (t) = \frac{\nu_0}{4} | A\u^n |^2, \quad
h(t) =  C_\eps  ( 1   + |\nabla m_\alpha|_4^8 +  |\nabla m_\alpha|_4^2  + C_{\eta_0}^8  +C_{\eta_0}^2),
\vspace{0.2cm}
\\
\displaystyle
g(t) = C  G(t)/(2 C_{\nu_0}) + 2 C |\eta^n_0|^2
+ C C_{\eta_0}^2 
+C  |m_\alpha|^2 
+ C |\f (t)|^2 
\end{array}
\end{equation}

We also define  some other useful expressions:

By observing that $|D(\u^n_0)|^2  \leq C_D| A^{1/2} u_0|$ 
and  using the previous definitions,
with $C$ being the positive constant appearing  in the first term of the right-hand side of (\ref{332Z}), 
the expression for the condition appearing in Lemma  \ref{lema110} can be estimated by
\[
C\int_0^{t}(2\varphi_0+2\int_0^{s}g(\tau)\,d\tau)^2+\int_0^{t}h(s)\,ds \leq F(t) ,
\]

\noindent
where
\[
F(t) = C\int_0^{t}(2 C_D| A^{1/2} u_0|+2\int_0^{s}g(\tau)\,d\tau)^2+\int_0^{t}h(s)\,ds .
\]

We also observe that the right-hand side in the conclusion of Lemma \ref{lema110} can be estimated by:
\[
2 (\varphi_0 + \int_0^t g(s) ds \leq \gamma (t),
\]

\noindent
with
\[
\gamma (t)  =  2 (C_D |D(\u_0)|^2 + \int^t_0 g(s) ds ) .
\]

\vspace{0.2cm}
Next, from an interpolation inequality and the estimates (\ref{334}) and (\ref{339}), we have that
\begin{equation}
\label{340}
\begin{array}{c}
|\nabla\xi^n(t)|_4 
 \leq  
C_6 |\nabla\xi^n(t)|^{1/4} |A_1\xi^n(t)|^{3/4}
\vspace{0.2cm}
\\
 \leq 
C_6 H_1 (t;  \|D(\u^n)\|_{L^\infty(0, t; L^2 (\Omega))})^{1/8}  H_3 (t; \|D(\u^n)\|_{L^\infty(0, t; L^2)} ; \|\nabla\xi^n\|_{L^\infty(0, t; L^4)} ; \|A\u^n\|_{L^2(0, t; L^2)})^{3/8}
\vspace{0.2cm}
\\
\leq 
H_4(t) ,
\end{array}
\end{equation}

\noindent
where we defined
\[
H_4(t) =  C_6 H^{1/8}_1 (t, \gamma(t)) \,
H^{3/8}_3 (  t, \gamma(t), 1, 4 \gamma (t)/\nu_0 ).
\]

\vspace{0.2cm}
By using these previous results, we can now choose $T^*$.
For this, we observe that $F(t)$ and $H_4(t)$ are continuous functions of $t$ and that $F(0) = H_4(0) =0$;
thus, there is  $T^*> 0$ satisfying the following inequalities:
\begin{equation}\label{FH}
F(T^*)    \leq \frac12\ln \frac{3}{2}.
\end{equation}
\begin{equation}\label{FH2}
H_4(T^*)  < 1 .
\end{equation}

Then, since $F(t)$ and $H_4(t)$ are  nondecreasing functions, (\ref{FH})  and (\ref{FH2})  imply that
\begin{equation}\label{343}
F(t) < \frac12\ln \frac{3}{2} \qquad {\rm and} \qquad H_4(t) < 1, \qquad
{\rm for} \quad t \in [0, T^*].
\end{equation}

\vspace{0.2cm}
With these observations, we will prove estimate (\ref{xi1}).

Indeed, since $|\nabla\xi^n(0)|=0$, by continuity, (\ref{xi1}) is true for small times. 
Assume by contradiction that there exists $T^n_1$, with $0 < T^n_1 < T^*$, such that   $|\nabla\xi^n(t)|_4 < 1$ for $0 \leq t < T^n_1$
and 
\begin{equation}
\label{HypothesisForContradiction1}
|\nabla\xi^n(T^n_1)|_4 = 1. 
\end{equation}

Thus, from (\ref{332Z}) and our choice  of $\eps$ in (\ref{epsilon}), we conclude that, for $t \in [0,T^n_1]$, $\u^n$ satisfies  the following differential inequality:
\begin{equation}
\label{vG}
\begin{array}{c}
\displaystyle
\frac{d}{dt}|D(\u^n)|^2 
  + \frac{\nu_0}{4} | A\u^n |^2
  \vspace{0.2cm}
  \\
  \displaystyle
  \leq
C |D(\u^n)|^6 
+ C_\eps  \left( |\nabla\xi^n|^8_4 + |\nabla m_\alpha|_4^8 +  |\nabla m_\alpha|_4^2  + C_{\eta_0}^8  +C_{\eta_0}^2 \right) |D(\u^n)|^2  
  \vspace{0.2cm}
  \\
  \displaystyle
+  C  G(t)/(2 C_{\nu_0}) + 2 C |\eta^n_0|^2
+ C C_{\eta_0}^2 
+C  |m_\alpha|^2 
+ C |\f|^2 
\end{array}
\end{equation}

By taking into account (\ref{FH}), we can to apply the Lemma \ref{lema110} to (\ref{vG}), 
with the expressions given in (\ref{Lemma10Correspondence}) to obtain that
\begin{equation}
\label{Dvn}
|D(\u^n(t)|^2 + \frac{\nu_0}4 \int^t_0 |A\u^n|^2 ds 
\leq 
2 ( |D(\u^n_0)|^2 +  \int_0^t g(s) ds)
=
\gamma (t)
\leq 
\gamma (T^*), \qquad {\rm for} \quad t \in [0, T^n_1].
\end{equation}

This last result, together with (\ref{334}), (\ref{339}) and since $|\nabla\xi^n(t)|_4 \leq1$ for $0 \leq t \leq T^n_1$, 
gives that
\[
\begin{array}{l}
|\nabla \xi^n (t)|^2 \leq H_1 (t, \gamma (t)) \leq H_1 (T^*, \gamma (T^*)) ,
\vspace{0.2cm}
\\
|A_1 \xi^n(t)|^2 \leq H_3 (t, \gamma(t), 1, 4 \gamma(t)/\nu_0)) \leq H_3 (T^*, \gamma(T^*), 1, 4 \gamma(T^*)/\nu_0)).
\end{array}
\]

\noindent
for $0 \leq t \leq T^n_1$.
Therefore, (\ref{340}), (\ref{FH2}) and  (\ref{343})  imply that $|\nabla \xi^n(t)|_4 \leq H_4 (t) < 1$ for $t \in [0, T^n_1]$. 
In particular $|\nabla\xi^n(T^n_1)|_4 < 1$,
in contradiction with (\ref{HypothesisForContradiction1}).
 This contradiction proves that  $T_1^n \geq T^*$ and thus the first four estimates in  (\ref{xi1}) hold.
By using these estimates in (\ref{parttxi}) and (\ref{336}), we obtain respectively the fifth and the sixth estimates in (\ref{xi1}).
\hfill $\Box$.

\subsection{Proof of Proposition \ref{localsol}}

The results in (\ref{xi1}) give higher order estimates for $(\u^n, \xi^n)$ and, as consequence, the required higher order estimates for  $(\u^n, \eta^n)$. 
They suffice to take limit when $n \rightarrow \infty$ along a suitable subsequence to obtain a local strong solution
with the required regularity (\ref{r21}), (\ref{r22}) on $[0,T^*]$.

As far as the initial conditions, for the velocity,  it suffices to show that
$\displaystyle\lim\sup_{t \rightarrow 0^+} |D(\u(t))| \leq |D(\u_0)|$ because we already know that $\u(t) \rightarrow \u_0$ strongly in $\mathbf{L}^2(\Omega)$. 
In fact, $\u$ is obviously a weak solution, and this convergence was proved in the previous section.

From (\ref{vG}), by using (\ref{Dvn}) and that $|\nabla \xi^n(t)|_4 \leq H_4 (t) < 1$ for $t \in [0, T^*]$,
we conclude that
\[
\begin{array}{c}
\displaystyle
|D(\u^n) (t)|^2 
  \leq
|D(\u^n_0)|^2 
+ \int_0^t \tilde{g} (s) \, ds ,
\end{array}
\]

\noindent
where $\tilde{g} (s) = C \gamma (s)^3 
+ C_\eps  \left( H_4 (s)^8  + |\nabla m_\alpha|_4^8 +  |\nabla m_\alpha|_4^2  + C_{\eta_0}^8  +C_{\eta_0}^2 \right) \gamma (s)  
+  C  G(t)/(2 C_{\nu_0}) + 2 C |\eta^n_0|^2
+ C C_{\eta_0}^2 
+C  |m_\alpha|^2 
+ C |\f (s)|^2 $.
By taking $n \rightarrow +\infty$ and then $t \rightarrow 0^+$ in this last inequality, we obtain the desired 
result. 

Also, since  $\eta(t) \rightarrow \eta_0$ in $L^2$ when $t \rightarrow 0^+$, 
by proceeding similarly as in the case of $\u$,  estimate (\ref{339}) implies that $|A_1\eta(t)|$
 remains bounded, the weak convergence $\eta(t) \rightarrow \eta_0$  in
$H^2 (\Omega)$ is obtained. This completes the proof of Proposition \ref{localsol}.
\hfill $\Box$

\vspace{0.2cm}
Therefore, once Proposition \ref{localsol} is proved, we guarantee the existence of local strong solutions as mentioned in Theorem \ref{localsolMainResults}.


\subsection{Proof of the uniqueness of strong evolution solutions}
\label{ProofOfTheUniquenessOfLocalEvolutionaryStrongSolutions}

Suppose that $(\u_1, m_1)$ and $(\u_2, m_2)$ are strong solutions with the same initial and boundary conditions, same external forces and same $\alpha$.
We define
\[
\v=\u_2-\u_1 \quad{and} \quad \eta=m_2-m_1  \quad{and} \quad q=q_2-q_1  .
\]

\noindent
We stress that here symbol $\eta$ has a different meaning than in the previous sections.

Then,  system (\ref{sistema2}) is rewritten in terms of $\v$, $\eta$ and $q$ as follows:
\begin{equation}
\label{sistema3Unicity}
  \left\{
  \begin{array}{l}
  \displaystyle\frac{\partial\v}{\partial t}-2\,{\rm div}
  (\nu(\eta+m_1)D(\v))-2\,{\rm div}
  (\nu(\eta+m_1)D(\u_1))
\\
  +2\,{\rm div}(\nu(m_1)D(\u_1))+
  \v\cdot\nabla \v+\v\cdot\nabla \u_1+
  \u_1\cdot\nabla \v+\nabla q
=-\eta\cdot \chi,
\\
  {\rm div}\, \v=0,
\\
  \displaystyle\frac{\partial\eta}{\partial t}-\theta\Delta\eta
  +\v\cdot\nabla\eta+\v\cdot\nabla
  m_1 +\u_1\cdot\nabla\eta+U
  \displaystyle\frac{\partial\eta}{\partial x_3}=0 
  \ {\rm in}\ (0,T)\times\Omega ,
  \end{array}
  \right.
  \end{equation}
  
  \noindent
  with boundary conditions
  \begin{equation}
  \label{bcPerturbationUnicity}
  \left\{
  \begin{array}{l}
  \v=0 \ {\rm on}\ (0,T)\times S,
  \\
  \v\cdot\textbf{n} =0 \ {\rm on}\ (0,T)\times \Gamma,
  \\
  \nu(\eta+m_1)[D(\v)\n-\n\cdot(D(\v)\n)\n]
  \\
  \displaystyle
  \hspace{2cm}
  =
  -  (\nu (\eta + m_1 ) - \nu (m_1) ) [D(\u_1)\n-\n\cdot(D(\u_1)\n)\n]
 ,
  \ {\rm on}\ (0,T)\times \Gamma,
  \\
  \theta\displaystyle\frac{\partial\eta}{\partial\n}-U\eta \, n_3 =0\ {\rm on}\ (0,T) ,
  \times\partial\Omega ,
  \end{array}
  \right.
  \end{equation}

 \noindent
and the following initial conditions
\begin{equation}
\label{ic3Unicity}
\v= 0, \quad \eta=0 \quad\mbox{at time} \; t = 0.
\end{equation}

\vspace{0.1cm}
In weak formulation, the previous problem becomes:
\begin{equation}
\label{12Unicity}
  \left\{
\begin{array}{l}
\displaystyle
 \left< \partial_t \v, \w\right> + 2(\nu(\eta+m_1)D(\v), D(\w)) + 2(\nu(\eta+m_1) D(\u_1), D(\w)) 
  \\\noalign{\medskip}
  - 2 (\nu(m_1) D(\u_1), D(\w)) + (\v\cdot\nabla\v,\w) +  (\v\cdot\nabla \u_1, \w) + (\u_1 \cdot\nabla \v, \w) =  
 \\
 - (\eta \cdot \chi, \w)
+ 2 \int_{\Gamma}  (\nu (\eta + m_1 ) - \nu (m_1) ) [D(\u_1)\n]_{\mathbf{t}} \cdot \w~dS ,
\vspace{0.1cm}
  \\
  \displaystyle
 \left< \partial_t \eta, \phi \right> + \theta (\nabla\eta, \nabla\phi) + (\v\cdot\nabla\eta,\phi) 
+ (\v\cdot\nabla m_1, \phi) +   (\u_1\cdot\nabla \eta, \phi) - U\left(\eta,
\frac{\partial\phi}{\partial x_3} \right) = 0,
\end{array}
\right.
\end{equation}
in $D'(0,T)$, for all $\w \in \J_0$ and $\phi \in B$.
Here, we denoted 
\[
[D(\u_1)\n]_{\mathbf{t}} = D(\u_1)\n-\n\cdot(D(\u_1)\n)\n .
\]

To prove uniqueness of strong solutions, it is enough to prove the exponential decay of the perturbations $(\v, \eta)$.

For this, we first take $\w = \v$ in the first equation of (\ref{12Unicity}),
using that, similarly as before, $\displaystyle \left< \partial_t \v, \v \right> =  \frac12 \frac{d}{dt} |\v|^2$,
to obtain that
\begin{eqnarray}
&& \frac12 \frac{d}{dt} |\v|^2 + 2 (\nu(\eta+m_1) D(\v), D(\v)) 
\nonumber
\\
&&
+ 2 \left( (\nu(\eta+m_1) - \nu(m_1)) D(\u_1), D(\v) \right) 
 + (\v \cdot\nabla \u_1, \v) + (\eta \cdot \chi, \v) 
\nonumber
\\
&& = 
 2 \int_{\Gamma}  (\nu (\eta + m_1 ) - \nu (m_1) ) [D(\u_1)\n]_{\mathbf{t}} \cdot \v~dS .
\label{vwnUnicity}
\end{eqnarray}

Next, by using \eqref{nus}, H\"older and Young inequalities and Sobolev embeddings, we obtain the following estimates, for the term in (\ref{vwnUnicity}):
\[
 2 (\nu(\eta+m_1) D(\v), D(\v)) \geq 2 \nu_0 |D(\v)|^2 ,
\]

\[
\begin{array}{c}
 2 | \left( (\nu(\eta+m_1) - \nu(m_1)) D(\u_1), D(\v) \right) | \leq 2 \nu^\prime_1 |\eta|_4 \, |D(\u_1)|_4 \ |D(\v)|
\vspace{0.2cm}
\\
\leq  2  C_1 \nu^\prime_1 \|\u_1\|_2 | \nabla \eta| \ |D(\v)|  
\vspace{0.2cm}
\\
\displaystyle
\leq
  C_{\nu_0} 4  C_1^2 {\nu^\prime_1}^2  \|\u_1\|_2^2 \, | \nabla \eta|^2 
+ \frac{1}{2} \nu_0 |D(\v)|^2 ,
\end{array}
\]

\[
\begin{array}{c}
|(\v \cdot\nabla \u_1, \v) | \leq |\v|_4 |\nabla \u_1| |\v |_4 \leq C  \| \u_1\|_1 |\v |_4^2 \leq C  \| \u_1\|_1 |\v |^{1/2}  | \nabla \v |^{3/2}
\vspace{0.2cm}
\\
\displaystyle
\leq C_{\nu_0} C  \| \u_1\|_1^4 |\v |^2 + \frac{1}{2} \nu_0 |D(\v)|^2 ,
\end{array}
\]

\[
| (\eta \cdot \chi, \v) | \leq |\eta| | \v | \leq  | \v |^2 + |\nabla \eta|^2 ,
\]

\[
\begin{array}{c}
2 \left|\displaystyle\int_{\Gamma} -  (\nu (\eta + m_1 ) - \nu (m_1) ) [D(\u_1)\n]_{\mathbf{t}} \cdot \v~dS \right| 
\leq
2 \nu_1^\prime \| \eta \|_{L^4 (\Gamma)} \| D (\u_1 ) \|_{L^2 (\Gamma)} \|\v \|_{L^4 (\Gamma)}
\vspace{0.2cm}
\\
\leq
2 C \nu_1^\prime | \nabla \eta | \, \|\u_1 \|_2 \,  | D(\v) |
\vspace{0.2cm}
\\
\displaystyle
\leq
4 C_{\nu_0} C^2 {\nu_1^\prime}^2 \,  \|\u_1 \|_2^2 \, | \nabla \eta |^2  + \frac{1}{2} \nu_0 | D(\v) |^2 .
\end{array}
\]

By plugging these last estimates in (\ref{vwnUnicity}) and simplifying, we obtain that
\begin{equation}
\label{UnicityPreparatoryv}
\frac{d}{dt} |\v|^2 +  \nu_0 |D(\v)|^2  
\leq  
\left( C_2{\nu^\prime_1}^2  \| \u_1\|_1^4   + 1 \right)  | \v |^2
+ \left(C_2 {\nu^\prime_1}^2  \|\u_1\|_2^2 +1 \right) |\nabla \eta|^2 .
\end{equation}

\vspace{0.2cm}
Next,  by taking $\phi = - A_1 \eta$ in the second equation in (\ref{12Unicity}), we obtain that
\begin{equation}
\label{EtaUnicity}
\begin{array}{c}
 \displaystyle
\frac{d}{dt} |\nabla \eta|^2
+ \theta |A_1 \eta |^2
\leq
C | \v\cdot\nabla\eta |^2
+ C |\v\cdot\nabla m_1 |^2
+  C | \u_1\cdot\nabla \eta |^2
+  U C | \nabla \eta|^2 .
\end{array}
\end{equation}

By using interpolation and Young inequality, we can estimate the terms in (\ref{EtaUnicity}) as follows.
\[
\begin{array}{l}
C | \v\cdot\nabla\eta |^2 \leq C  | \v|_4^2 |\nabla\eta |_4^2 \leq C | \v|^{1/2}|\nabla \v|^{3/2}  |\nabla\eta |_4^2 
\vspace{0.2cm}
\\
\hspace{1.55cm}
\leq C  \|\eta \|_2^8 | \v|^2+  | \nabla \v|^2 
\vspace{0.2cm}
\\
\hspace{1.55cm}
\leq C  ( \|m_1 \|_2 + \|m_2 \|_2)^8 | \v|^2+  | D(\v) |^2  ,
\end{array}
\]

\[
C  |\v|_4^2 |\nabla m_1 |_4^2 \leq C   \| m_1 \|_2^2 |D(\v)|^2 ,
\]

\[
C | \u_1\cdot\nabla \eta |^2 \leq C  | \u_1|_\infty^2 |\nabla \eta |^2 \leq  C \| \u_1 \|_2^2 |\nabla \eta |^2 .
\]

By pluggin these last estimates in (\ref{EtaUnicity}) and simplifying, we obtain that
\begin{equation}
\label{UnicityPreparatoryEta}
\begin{array}{c}
 \displaystyle
\frac{d}{dt} |\nabla \eta|^2
+ \theta |A_1 \eta |^2
\leq
  C  \left( \|m_1 \|_2 + \|m_2 \|_2 \right)^8 | \v|^2
\vspace{0.2cm}
\\
+ C  \left( \| \u_1 \|_2^2 +  U \right) | \nabla \eta|^2 
\vspace{0.2cm}
\\
+ \left( C   \| m_1 \|_2^2  + 1 \right) | D(\v) |^2 .
\end{array}
\end{equation}

By denoting the constant
\[
K= \left( C   \| m_1 \|_{L^\infty(0,T^*, H^2)}^2  + 1 \right) ,
\]

\noindent
multiplying (\ref{UnicityPreparatoryv}) by  $2 K/\nu_0$
and adding to (\ref{UnicityPreparatoryEta}), we get that
\begin{equation}
\begin{array}{c}
 \displaystyle
\frac{d}{dt} \left[ (2 K/\nu_0) |\v|^2   + |\nabla \eta|^2 \right]
+ \theta |A_1 \eta |^2
+  K  |D(\v)|^2  
\vspace{0.2cm}
\\
\leq  
N(t) \left[  (2 K/\nu_0)  | \v|^2 + | \nabla \eta|^2 \right] ,
\end{array}
\end{equation}

\noindent
where
\[
\begin{array}{l}
N(t) = 
 C_2{\nu^\prime_1}^2  \| \u_1\|_1^4   + 1 + C  \left( \|m_1 \|_2 + \|m_2 \|_2 \right)^8   (2 K/\nu_0)^{-1}
\vspace{0.2cm}
\\
\hspace{1.2cm}
+ 
(2 K/\nu_0) (C_2 {\nu^\prime_1}^2  \|\u_1\|_2^2 +1 )  + C  ( \| \u_1 \|_2^2 +  U ) .
\end{array}
\]

Since $N(\cdot)$ is an integrable function and $\v(0) = 0$ and $\eta(0) = 0$, using the Gronwall lemma, we obtain the
desired result.
\hfill $\Box$

\vspace{0.3cm}
Thus, the results of this section prove Theorem  \ref{localsolMainResults}.

\section{Existence of global strong evolution solutions}
\label{GlobalEvolutionaryStrongSolutions}

Imposing smallness conditions on the problem data, it is possible to prove that the local strong evolution solutions obtained in the previous sections are global in time.

\begin{theorem} 
\label{Teorema 4.1.MainResults} 
Let $m_\alpha$ be a solution of Lemma \ref{AuxiliaryProblemLemma}.
Suppose that the assumptions of Theorem \ref{localsolMainResults} hold,
with $\f \in L^\infty (0,\infty; \X)$ and $m_0 - m_\alpha \in Dom(A_1)  \subset Y \cap H^2 (\Omega)$.
Additionally, assume that the norms  
$|D(\u_0)|$, $|\u_0|_4$, $|\Delta (m_0 - m_\alpha )|$, $|\nabla (m_0 - m_\alpha )|_4$, 
$|\nabla (m_0 - m_\alpha )|$, $|m_0 - m_\alpha |$,  $U$ and $\| \f \|_{L^\infty (0,\infty; \X)} $
are small enough.
Then the solution $(\u,m)$ given in Theorem \ref{localsolMainResults}  exists globally in time. 
Moreover, there exist positive constants $M_1$ and $M_2$, such that
\begin{equation}
\sup_{t\geq 0} |\Delta m (t)|^2 \leq M_1, \qquad \sup_{t\geq
0} \{|D(\u(t))| + |\partial_t m (t)|\} \leq M_2 .
\end{equation}

\end{theorem}

\begin{remark}
The results presented in Theorem \ref{Teorema 4.1.MainResults}  are proved by using spectral Galerkin approximations.
Along the proof,  estimates similar to the stated ones are obtained holding uniformly for the approximate solutions;
these estimates can then be used to obtain the associated error estimates.
\end{remark}

\subsection{Proof of existence of global strong evolution solutions}

From the change of variables (\ref{ChangeOFVariablesConcentrationOther}), to prove Theorem \ref{Teorema 4.1.MainResults}, it is enough to prove the following result.

\begin{proposition} 
\label{Teorema 4.1.} 
Suppose that the assumptions of Theorem \ref{localsolMainResults} hold
with $\f \in L^\infty (0,\infty; \X)$ and $\eta_0 \in Dom(A_1) \subset Y \cap H^2 (\Omega)$.
Additionally, that the norms
$|D(\u_0)|$, $|\u_0|_4$, $|\Delta \eta_0|$, $|\nabla \eta_0|_4$, $|\nabla \eta_0|$, $|\eta_0|$, $U$ and $\| \f \|_{L^\infty (0,\infty; \X)} $ are small enough.
Then the solution $(\u,\eta)$ given in Proposition \ref{localsol}
exists globally. Moreover, there exist positive constants $\widetilde{M}_1$
and $\widetilde{M}_2$, such that
\begin{equation}
\sup_{t\geq 0} |\Delta\eta (t)|^2 \leq \widetilde{M}_1, \qquad \sup_{t\geq
0} \{|D(\v(t))| + |\partial_t \eta (t)|\} \leq \widetilde{M}_2 .
\end{equation}

\noindent
Analogous estimates hold uniformly for the spectral Galerkin approximations for large $n$.
\end{proposition}

\noindent 
{\bf Proof:} 
By using the change of variables (\ref{FurtherChangeOfVariableConcentration}),
estimate (\ref{332Preparatory}) implies that
\[
\begin{array}{c}
\displaystyle
\frac{d}{dt}|D(\u^n)|^2 
  + \frac{\nu_0}{2} | A\u^n |^2
  \vspace{0.2cm}
  \\
  \displaystyle
  \leq
  \displaystyle
 C |D(\u^n)|^6 
+ C \left( ( |\nabla \eta^n|_4^2 +  |\nabla m_\alpha|_4^2 )^4 +
 ( |\nabla \eta^n |_4 + |\nabla  m_\alpha |_4 )^2  \right) |D(\u^n)|^2  
  \vspace{0.2cm}
  \\
  \displaystyle
+ \widetilde{C}  | \eta^n |^2  + C  | m_\alpha |^2 
+ C |\f|^2 
+  \left(  |\nabla \eta^n |_4 + | \nabla m_\alpha |_4 \right)  |A\u^n|^2 .
\end{array}
\]

To simplify the notation, we denote
\begin{equation}
\label{DefinitionPi}
\Pi_n (t) =  |A_1 \eta^n (t) |^2 + |\nabla  m_\alpha |_4^2 .
\end{equation}

\noindent
Since from (\ref{EstimateImportant}) we have that
\[
 |\nabla\eta^n|_4^2 \leq \hat{C} |A_1 \eta|^2, 
\]

\noindent
with $ \hat{C} > 0$ independent of $\eta^n$, we obtain that
\[
|\nabla \eta^n |_4^2 + |\nabla  m_\alpha |_4^2 \leq \max \{ 1, \hat{C} \} \, \Pi_n (t).
\]

So, our last estimate implies that
\begin{equation}
\label{D(un)GS}
\begin{array}{l}
\displaystyle
\frac{d}{dt}|D(\u^n)|^2 
  + \frac{\nu_0}{2} | A\u^n |^2
  \leq
 C \left(
|D(\u^n)|^6 
+  \left( \Pi_n (t)^4 + \Pi_n (t)  \right)  |D(\u^n)|^2  
+\Pi_n (t)^{1/2} |A\u^n|^2 
\right)
  \vspace{0.2cm}
  \\
  \displaystyle
\hspace{4.4cm}
+  \widetilde{C}  | \eta^n |^2  
+ C  | m_\alpha |^2 
+ C |\f|^2 .
\end{array}
\end{equation}

\vspace{0.1cm}
However, we need further estimates to prove the proposition.

For this, we take $\phi = \eta^n (\cdot, t)$ in the second equation of (\ref{ApproximateEquationsOthers});
after some standard computations, we obtain
\[
\begin{array}{l}
\displaystyle
\frac{1}{2} \frac{d}{d t} |\eta^n |^2   + \theta |\nabla \eta^n |^2 
 \leq   
  |\u^n |_4^2 |\nabla m_\alpha |_4^2
+ U | \eta^n | \, | \nabla \eta^n |
   + U |m_\alpha| \, | \nabla \eta^n |
\vspace{0.2cm}
\\
\hspace{3.4cm}
\displaystyle
 \leq   
 |\nabla m_\alpha |_4^2   |D (\u^n )|^2  
+ 2 U C_p \, | \nabla \eta^n |^2
+ \frac{\theta}{4} |\nabla \eta^n |^2 .
\end{array}
\]

Now, we take $U$ small enough such that
\begin{equation}
\label{SmallnessConditionOnUGS1}
2 U C_P \leq \frac{\theta}{4},
\end{equation}

\noindent
and observe that, from (\ref{DefinitionPi}),
\[
  |\nabla m_\alpha |_4^2   |D (\u^n )|^2   \leq \Pi_n (t) | A \u^n |^2 .
\]

\noindent
By using this in the last estimate, we conclude that
\begin{equation}
\label{EtaGS}
 \frac{d}{d t} |\eta^n |^2   + \theta |\nabla \eta^n |^2 
 \leq   
 \Pi_n (t) | A \u^n |^2  .
\end{equation}

Since 
\[
\widetilde{K} |\eta^n|^2 \leq |\nabla \eta^n|^2, 
\]

\noindent
with $\widetilde{K} > 0$ independent of $\eta^n$,
by adding (\ref{D(un)GS}) and (\ref{EtaGS}) multiplied by 
\[
\hat{K} = \frac{ 2 \widetilde{C}} {\theta \widetilde{K} } ,
\]

\noindent
after some simplification, we obtain that
\begin{equation}
\label{D(un)EtanGSPreparatory}
\begin{array}{c}
\displaystyle
\frac{d}{dt}\left( |D(\u^n)|^2  + \hat{K} |\eta^n |^2  \right)
  + \frac{\nu_0}{2} | A\u^n |^2
 + \frac{\hat{K}  \theta }{ 2} |\nabla \eta^n |^2 
  \vspace{0.2cm}
  \\
  \displaystyle
  \leq
 C  |D(\u^n)|^6 
+ C \left( \Pi_n (t)^4 + \Pi_n (t)  \right)  |D(\u^n)|^2  
+C \Pi_n (t)^{1/2} |A\u^n|^2 
  \vspace{0.2cm}
  \\
  \displaystyle
+ C |\f|^2 
+ C  |m_\alpha |^2 .
\end{array}
\end{equation}

Next, we differentiate the second equation in (\ref{ApproximateEquationsOthers})  with respect to $t$ and take $\phi = \partial_t \eta^n (\cdot,t)$ 
to obtain, after some standard computations, that
\[
\begin{array}{c}
\displaystyle
\frac{1}{2} \frac{d}{d t} | \partial_t \eta^n|^2
+ \theta |\nabla  \partial_t  \eta ^n|^2
=
  - ( \partial_t  \u^n\cdot\nabla m_\alpha, \partial_t \eta^n  )
 - (\partial_t  \u^n \cdot\nabla \eta^n , \partial_t \eta^n ) 
 + U\left( \partial_t  \eta^n  , \frac{\partial}{\partial x_3}  \partial_t \eta^n \right) .
\end{array}
\]

The terms in the right-hand side can be estimated as
\[
 | ( \partial_t  \u^n\cdot\nabla m_\alpha, \partial_t \eta^n  ) | 
\leq
|\partial_t  \u^n | \, |\nabla m_\alpha |_4 \, | \partial_t \eta^n ) | _4
 \leq 
C_\theta |\nabla m_\alpha |_4^2   |\partial_t  \u^n |^2  + \frac{\theta}{8}  |\nabla \partial_t \eta^n ) |^2 ,
\]

\[
|(\partial_t  \u^n \cdot\nabla \eta^n , \partial_t \eta^n ) | 
\leq 
|\partial_t  \u^n | \, |\nabla \eta^n |_4 \, | \partial_t \eta^n ) | _4
 \leq 
C_\theta |\nabla \eta^n |_4^2   |\partial_t  \u^n |^2  + \frac{\theta}{8}  |\nabla \partial_t \eta^n ) |^2 ,
\]

\[
| U\left( \partial_t  \eta^n  , \frac{\partial \partial_t \eta^n}{\partial x_3} \right) | 
\leq
 U C_P | \partial_t  \eta^n | | \nabla \partial_t \eta^n|
\leq 
 U \frac{ C_P}{4} | \partial_t  \eta^n |^2  + U C_P | \nabla \partial_t \eta^n|^2 .
\]

By using these last results and the smallness condition on $U$, we then obtain that
\[
\begin{array}{c}
\displaystyle
 \frac{d}{d t} | \partial_t \eta^n|^2
+ \theta |\nabla  \partial_t  \eta ^n|^2
\leq
2 C_\theta  |\nabla \eta^n |_4^2   |\partial_t  \u^n |^2 
 + \left( |\nabla m_\alpha |_4^2 +  U \frac{ C_P}{2} \right) | \partial_t  \eta^n |^2  .
\end{array}
\]

Since there is $\tilde{C} >0$ independent of $\eta^n$ such that $ \tilde{C} |\partial_t  \eta^n |^2  \leq |\nabla  \partial_t  \eta ^n|^2$,
recalling (\ref{MalphaEstimates}) and taking $U$ is small enough such that, besides (\ref{SmallnessConditionOnUGS1}), the following also holds
\begin{equation}
\label{SmallnessConditionOnUGS2}
\left( |\nabla m_\alpha |_4^2 +  U \frac{ C_P}{2} \right)  \leq \frac{\theta  \tilde{C}}{2} ,
\end{equation}

\noindent
from the last estimate we obtain that
\begin{equation}
\label{PartialtEtaGS}
\begin{array}{c}
\displaystyle
 \frac{d}{d t} | \partial_t \eta^n|^2
+ \frac{\theta}{2} |\nabla  \partial_t  \eta ^n|^2
\leq
2 C_\theta |\nabla \eta^n |_4^2   |\partial_t  \u^n |^2 .
\end{array}
\end{equation}

Next, we put $\w = \partial_t  \u^n (\cdot, t)$ in the first equation of (\ref{ApproximateEquationsOthers}) .
Since we can estimate
\[
\begin{array}{l}
|\u^n\cdot\nabla\u^n |^2
 \leq |\u^n |_6^2 |\nabla\u^n |_3^2 
\leq C |D(\u^n )|^2 |D(\u^n )| |A\u^n |
\vspace{0.2cm}
\\
\hspace{2cm}
\leq C  |D(\u^n) |^3 |A \u^n |
\leq C  |D(\u^n) |^6 + |A \u^n |^2 ,
\end{array}
\]

\noindent
after integration by parts and some rather standard computations, we obtain that
\begin{equation}
\label{PartialtuGS}
\begin{array}{l}
\displaystyle
 |\partial_t \u^n |^2
\leq
C \left( 
|A \u^n |^2
  + |m_\alpha |^2
+ |\f |^2
\right)
+ C |D(\u^n) |^6
+ |A \u^n |^2
 +\frac{1}{2} |\eta^n|^2
\vspace{0.2cm}
\\
\hspace{1.3cm}
\displaystyle
\leq
 C |D(\u^n) |^6
+ C |A \u^n |^2
+C  |m_\alpha |^2
+ C |\f |^2
 +\frac{1}{2} |\eta^n|^2 .
\end{array}
\end{equation}

By plugging (\ref{PartialtuGS}) in (\ref{PartialtEtaGS}), recalling the definition (\ref{DefinitionPi}) of $\Pi (t)$, we obtain that
\[
\begin{array}{c}
\displaystyle
 \frac{d}{d t} | \partial_t \eta^n|^2
+ \frac{\theta}{2} |\nabla  \partial_t  \eta ^n|^2
\leq
2 C_\theta \Pi (t) 
\left( 
 C |D(\u^n) |^6
+ C |A \u^n |^2
+C  |m_\alpha |^2
+ C |\f |^2
 +\frac{1}{2} |\eta^n|^2 
 \right) .
\end{array}
\]

By adding this last inequality to (\ref{D(un)EtanGSPreparatory}), after some simplifications and rearranging, we get that
\begin{equation}
\label{DifferentialInequalityPreparatory}
\begin{array}{c}
\displaystyle
\frac{d}{dt}\left( |D(\u^n)|^2  + \hat{K} |\eta^n |^2  +  | \partial_t \eta^n|^2 \right)
  + \frac{\nu_0}{2} | A\u^n |^2
 + \frac{\hat{K}  \theta }{ 2} |\nabla \eta^n |^2 
+ \frac{\theta}{2} |\nabla  \partial_t  \eta ^n|^2
  \vspace{0.2cm}
  \\
  \displaystyle
  \leq
 \hat{C}_1 ( 1 + \Pi_n (t) )  |D(\u^n)|^6 
+ \hat{C}_2  \left( \Pi_n (t)  + \Pi_n (t)^4  \right)  |D(\u^n)|^2  
+  \hat{C}_3  \Pi_n (t)  |\eta^n|^2
  \vspace{0.2cm}
  \\
  \displaystyle
+  \hat{C}_4  \Pi_n (t)^{1/2} |A\u^n|^2 
+  \hat{C}_5  (1 + \Pi_n (t) )|\f|^2 
+  \hat{C}_6  \left ( 1  + \Pi_n (t) \right) |m_\alpha |^2 ,
\end{array}
\end{equation}

\noindent
for certain positive constants $\hat{C}_i$, $i=1, \ldots, 6$.

Next, besides the definition of $\Pi (t)$ given in (\ref{DefinitionPi}), we also denote
\begin{equation}
\label{DefinitioVarphiPsi}
\left\{
\begin{array}{l}
\displaystyle
\varphi_n (t) = |D(\u^n)|^2  + \hat{K} |\eta^n |^2  +  | \partial_t \eta^n|^2  ,
\vspace{0.2cm}
\\
\displaystyle
\psi_n (t) = \frac{\nu_0}{4} | A\u^n |^2
 + \frac{\hat{K}  \theta }{ 2} |\nabla \eta^n |^2 
+ \frac{\theta}{2} |\nabla  \partial_t  \eta ^n|^2 .
\end{array}
\right.
\end{equation}

\noindent
Then,  estimate (\ref{DifferentialInequalityPreparatory}) implies that
\begin{equation}
\label{DifferentialInequality1}
\begin{array}{l}
\displaystyle
\frac{d \varphi_n (t)}{dt} + \psi_n (t) 
+ \frac{\nu_0}{4} | A\u^n (t) |^2 
  \leq
 D_1 ( 1 + \Pi_n (t) ) \varphi (t)^3 
+ D_2 ( \, \Pi_n (t)  + \Pi_n  (t)^4 )   \,\varphi (t)
  \vspace{0.2cm}
  \\
  \displaystyle
\hspace{2.5cm}
+  \hat{C}_4  \Pi_n (t)^{1/2} |A\u^n|^2 
+  \hat{C}_5  (1 + \Pi_n (t) )|\f|^2 
+  \hat{C}_6  \left ( 1  + \Pi_n (t) \right) |m_\alpha |^2 ,
\end{array}
\end{equation}

\noindent
for some positive constants $D_1$ and $D_2$.

Moreover, since by Sobolev embeddings, there is a positive constant $D_3$, independent of $n$, such that
\begin{equation}
D_3 \varphi_n (t) \leq \psi_n (t) ,
\end{equation}

\noindent
by using this in (\ref{DifferentialInequality1}), we get that
\begin{equation}
\label{DifferentialInequality2}
\begin{array}{c}
\displaystyle
\frac{d \varphi_n (t) }{dt}
+ \frac{\nu_0}{4} | A\u^n (t) |^2 
  \leq
 D_1 (1 + \Pi_n (t) ) \varphi_n (t)^3 
+ \left(  D_2 ( \, \Pi_n (t)  + \Pi_n (t)^4 )  - D_3 \right) \varphi_n (t)
  \vspace{0.2cm}
  \\
  \displaystyle
+  \hat{C}_4  \Pi_n (t)^{1/2} |A\u^n|^2 
+  \hat{C}_5  (1 + \Pi_n (t) )|\f|^2 
+  \hat{C}_6  \left ( 1  + \Pi_n  (t) \right) |m_\alpha |^2 .
\end{array}
\end{equation}

\vspace{0.2cm}
Next, by putting $\phi = - A_1 \eta^n (\cdot, t) $ in the second equation of  (\ref{ApproximateEquationsOthers}), 
proceeding similarly as we did in (\ref{StrongEstimateConcentration1}), using integration by parts, (\ref{SmallnessConditionOnUGS1}) and some rather standard computations,
we have that
\begin{equation}
 \label{A1EtaEstimateGS}
\begin{array}{l}
\displaystyle
 |  A_1 \eta^n |^2
\leq
\frac{2}{ \theta}
C_\theta
\left(
  | \partial_t \eta^n |^2
 + | \u^n |_4^2 | \nabla \eta^n |_4^2
+ |\u^n|_4^2 |\nabla m_\alpha |_4^2
 + U | \nabla  \eta^n |^2
\right)
\vspace{0.2cm}
\\
\hspace{1.4cm}
\displaystyle
\leq
\frac{2}{ \theta}
C_\theta C
\left( \,
  | \partial_t \eta^n |^2
 +   | D(\u^n ) |^2 |  |  A_1 \eta^n |^2 
+ | D(\u^n ) |^2  |\nabla m_\alpha |_4^2 
+ U |A_1 \eta^n|^2
\, \right)
\vspace{0.2cm}
\\
\hspace{1.4cm}
\displaystyle
\leq
D_4
 \varphi_n (t) \left(
 1
 + |\nabla m_\alpha |_4^2 
\right)
+ D_4 \left( \, \varphi_n (t)  +   U  \, \right) | A_1 \eta^n |^2 ,

\end{array}
\end{equation}

\noindent
where we used the definition of $\varphi (t)$ given in (\ref{DefinitioVarphiPsi}).

\vspace{0.2cm}
Next, we fix $\beta > 0$  small enough such that
\begin{equation}
\label{ConditionsOnBeta}
\begin{array}{l}
\displaystyle
\beta \leq 1, \quad
 \hat{C}_4  \beta^{1/2} \leq \frac{\nu_0}{4} \quad \mbox{and} \quad D_2 ( \beta + \beta^4 )  \leq  \frac{ D_3}{2} ,
\end{array}
\end{equation}

\noindent
where $\hat{C}_4$, $D_2$ and $D_3$ are the positive constants independent of $n$ appearing in (\ref{DifferentialInequality2}).

\vspace{0.2cm}
Now, we consider the polynomial
\begin{equation}
\label{Polynomialp(z)}
p(z) = 2 D_1 z^3   - \frac{D_3}{2} z + D_5,
\end{equation}

\noindent
with
\[
D_5 = 2 \hat{C}_5 \|\f \|_{L^2 (0, +\infty; \L^2)}^2  + 2 \hat{C}_6 |m_\alpha |^2 ,
\]

\noindent
and $D_1$, $D_3$, $ \hat{C}_5$ and $ \hat{C}_6$ the positive constants in (\ref{DifferentialInequality2}),
and analyze its roots as functions of (small) $D_5$.

We observe that, when $D_5 = 0$, the polynomial $p(z)$ has three simple real roots: $z_1 (0) = - \frac{1}{2} \sqrt{\frac{D_3}{D_1} }$, $z_2 (0)= 0$ and $z_3 (0) = \frac{1}{2} \sqrt{\frac{D_3}{D_1} }$.
Since the addition of  $D_5$ just translates the graph of $p(z)$ in the vertical direction,
for small $D_5 > 0$, $p(z)$ still has three simple real roots $z_1 (D_5) < 0$, $z_2 (D_5) > 0$, $z_3 (D_5) > 0$ and  $ z_2 (D_5) \rightarrow 0+$ as $D_5 \rightarrow 0+$.

The root $z_2 (D_5) $, as we will see in the arguments that follows, plays an important role in the derivation of estimates that will prove the proposition.

Since $\beta > 0$ is already fixed satisfying the conditions stated in (\ref{ConditionsOnBeta}),
we can assume that $ \|\f \|_{L^2 (0, +\infty; \L^2)}$ and $U$  are small enough, and so also $|\nabla m_\alpha |_4$ is small, such that
\begin{equation}
\label{FurtherConditions}
\begin{array}{l}
\displaystyle
D_4 z_2 (D_5)
\left(
  1
 +  |\nabla m_\alpha |_4^2 
\right)
+ D_4 \left( \, z_2  (D_5) +   U  \, \right)  \beta
\leq 
\frac{\beta}{4}
,
\vspace{0.2cm}
\\
\displaystyle
|\nabla m_\alpha|_4^2  \leq \frac{\beta}{8} .
\end{array}
\end{equation}

\noindent
where $D_4$ is the positive constant also independent of $n$ appearing in (\ref{A1EtaEstimateGS}).

\vspace{0.2cm}
Next, we consider the smallness conditions to be imposed on the initial conditions.

For this, we first observe that by using $\partial_t \eta^n (\cdot, t)$ in the second equation of (\ref{ApproximateEquationsOthers}),
by doing using integration by parts and standard computations, we obtain that
\[
| \partial_t \eta^n |^2
\leq 4 \left(
   \theta^2 |\Delta \eta^n |^2
  + |\u^n |_4^2 |\nabla \eta^n|_4^2
  +  |\u^n|_4^2 |\nabla m_\alpha |_4^2
 + U^2 | \nabla \eta^n |^2
\right) .
\]

By computing at time $t=0$, we then get that
\[
| \partial_t \eta_0^n |^2
\leq 4 \left(
   \theta^2 |\Delta \eta_0^n |^2
  + |\u_0^n |_4^2 |\nabla \eta_0^n|_4^2
  +  |\u_0^n|_4^2 |\nabla m_\alpha |_4^2
 + U^2 | \nabla \eta_0^n |^2
\right) .
\]

By proceeding similarly with the second equation in (\ref{WeakFormulationEquationsOthers}), we obtain that
\[
| \partial_t \eta_0 |^2
\leq 4 \left(
   \theta^2 |\Delta \eta_0 |^2
  + |\u_0^n |_4^2 |\nabla \eta_0|_4^2
  +  |\u_0|_4^2 |\nabla m_\alpha |_4^2
 + U^2 | \nabla \eta_0 |^2
\right) .
\]

Then, by taking the initial conditions $\u_0$ and $\eta_0$ small enough as in the statement of the Proposition \ref{Teorema 4.1.}, we can assure that
\begin{equation}
\label{ConditionsOnIC}
\begin{array}{l}
\displaystyle
| \partial_t \eta_0 |^2
\leq 4 \left(
   \theta^2 |\Delta \eta_0 |^2
  + |\u_0^n |_4^2 |\nabla \eta_0|_4^2
  +  |\u_0|_4^2 |\nabla m_\alpha |_4^2
 + U^2 | \nabla \eta_0 |^2
\right)
\leq
\frac{z_2 (D_5)}{8} ,
\vspace{0.2cm}
\\
\displaystyle
|D(\u_0)|^2  + \hat{K} |\eta_0 |^2  \leq \frac{z_2 (D_5)}{8} .
\end{array}
\end{equation}

Since $\u_0^n$ and $\eta_0^n$ converges respectively in $\H^1(\Omega)$ and $H^2(\Omega)$ to $\u_0$ and $\eta_0$, see (\ref{InitialDataConvergences1}) and (\ref{InitialDataConvergences2}),
conditions
(\ref{ConditionsOnIC}) implies that, for $n$ large enough,
\begin{equation}
\label{ConditionOnvarph0}
\begin{array}{l}
\displaystyle
\varphi_n (0) = ( |D(\u_0^n)|^2  + \hat{K} |\eta_0^n |^2  ) +  | \partial_t \eta_0^n|^2
\vspace{0.2cm}
\\
\hspace{1.2cm}
\displaystyle
\leq  ( |D(\u_0^n)|^2  + \hat{K} |\eta_0^n |^2  ) 
\vspace{0.2cm}
\\
\hspace{1.5cm}
\displaystyle
+   4 \left(
   \theta^2 |\Delta \eta_0^n |^2
  + |\u_0^n |_4^2 |\nabla \eta_0^n|_4^2
  +  |\u_0^n|_4^2 |\nabla m_\alpha |_4^2
 + U^2 | \nabla \eta_0^n |^2
\right)
\vspace{0.2cm}
\\
\hspace{1.2cm}
\displaystyle
  \leq 2 \frac{z_2 (D_5)}{8}+ 2 \frac{z_2 (D_5)}{8}
  \leq \frac{z_2 (D_5)}{2}.
\end{array}
\end{equation}

Additionally, we assume that $\eta_0$ is such that
\[
|A_1 \eta_0 |^2  \leq \frac{\beta}{10} ,
\]

\noindent
which, together with the second condition in (\ref{FurtherConditions}), for large enough $n$, 
see  (\ref{InitialDataConvergences2}),
implies that 
\begin{equation}
\label{ConditionPin(0)}
\Pi_n (0) =  |A_1 \eta_0^n |^2 + |\nabla  m_\alpha |_4^2 \leq \frac{\beta}{8} + \frac{\beta}{8} = \frac{\beta}{4}.
\end{equation}

\vspace{0.2cm}
Now we are ready to prove the existence of global in time strong evolution solutions.

For this,  we observe that,  from (\ref{ConditionPin(0)}) and continuity, we have that $\Pi_n (t) < \beta$ at least for small times $t$.
We will show that, under the stated conditions, we have in fact that $\Pi_n (t) < \beta$ for all $t$.

To show this, suppose by contradiction that this is not true and so there is $ 0< T_2^n < +\infty$ such that
\begin{equation}
\label{DefinitionT2n}
\Pi_n (t) < \beta \quad \mbox{for} \; t \in [0,T_2^n) \quad \mbox{and} \quad  \Pi_n (T_2^n) = \beta.
\end{equation}

Then, from (\ref{DifferentialInequality2}) and due to the third condition in (\ref{ConditionsOnBeta}), for $t \in [0,T_2^n]$,   we obtain that
\[
\begin{array}{c}
\displaystyle
\frac{d \varphi_n (t) }{dt}
+ \frac{\nu_0}{4} | A\u^n (t) |^2 

  \leq
 D_1 ( 1+ \beta) \varphi_n (t)^3 
  - \frac{D_3}{2}  \varphi_n (t)
  \vspace{0.2cm}
  \\
  \displaystyle
+  \hat{C}_4  \beta^{1/2} |A\u^n (t)|^2 
+  \hat{C}_5  (1 + \beta )|\f (t)|^2 
+  \hat{C}_6  \left ( 1 + \beta \right) |m_\alpha |^2 ,
\end{array}
\]

\noindent
which, owing to the first and second condition in (\ref{ConditionsOnBeta}), implies that
\begin{equation}
\begin{array}{c}
\displaystyle
\frac{d \varphi_n (t)}{dt} 
  \leq
 2 D_1\varphi_n (t)^3 
  - \frac{D_3}{2}  \varphi_n (t)
+ 2 \hat{C}_5 \|\f \|_{L^2 (0, +\infty; \L^2)}^2 
+ 2 \hat{C}_6 |m_\alpha |^2
= p (\varphi_n (t)) ,
\end{array}
\end{equation}

\noindent
for $t \in [0,T_2^n]$, and
where $p (\cdot)$ is exactly the polynomial defined in (\ref{Polynomialp(z)}).

Thus, we can use results about differential inequalities, see Hale \cite{hale}, to conclude that, 
\[
\varphi_n (t) \leq z(t), \quad \mbox{for} \; t \in [0,T_2^n],
\]

\noindent
for any $z(t)$ satisfying
\begin{equation}
\label{DifferentialEquality}
\begin{array}{l}
\displaystyle
\frac{d z (t)}{dt} 
=
p(z (t)) ,
\vspace{0.2cm}
\\
z(0) \geq  \varphi_n (0).
\end{array}
\end{equation}

By taking in particular the constant solution $z(t) \equiv  z_2 (D_5)$ of (\ref{DifferentialEquality}), from 
(\ref{ConditionOnvarph0}), we obtain that
\[
\varphi_n (t) \leq z_2 (D_5), \quad \mbox{for} \; t \in [0,T_2^n] .
\]

\vspace{0.2cm}
By using this last result, from (\ref{A1EtaEstimateGS}) and conditions (\ref{FurtherConditions}), for $t \in [0,T_2^n]$, we get that
\begin{equation}
\begin{array}{l}
\displaystyle
 |  A_1 \eta^n (t) |^2
\leq
D_4  \varphi_n (t)
\left(
 1+  |\nabla m_\alpha |_4^2 
\right)
+ D_4 \left( \, \varphi_n (t)  +   U  \, \right)  \beta
\vspace{0.2cm}
\\
\hspace{1.35cm}
\displaystyle
\leq
D_4 z_2 (D_5)
\left(
  1
 +  |\nabla m_\alpha |_4^2 
\right)
+ D_4 \left( \, z_2  (D_5) +   U  \, \right)  \beta
\leq 
\frac{\beta}{4} .
\end{array}
\end{equation}

From this, and again  (\ref{FurtherConditions}), we obtain in particular that
\[
\Pi_n (T_2^n) =  |  A_1 \eta^n (T_2^n) |^2 + |\nabla m_\alpha|_4^2  \leq \frac{\beta}{4} + \frac{\beta}{8} < \beta ,
\]

\noindent
in contradiction with (\ref{DefinitionT2n}).

Thus, under the stated conditions, the approximate solutions are defined on $[0,+\infty)$, 
and the required estimates hold uniformly with respect to large $n$.

Passing to the limit in standard way we then obtain a global strong solution to the problem satisfying the stated estimates.
\hfill $\Box$

\section{Exponential $L^2$-stability}
\label{ExponentialL2StabilityOfSmallStrongStationarySolutions}

Since in this section we are interested in the relationship between stationary solutions and solutions evolving in
time, here we consider only {\bf external forces $\f$ independent of time}.
Moreover, we analyze the evolution of disturbances around a strong stationary solution.

\subsection{Definition and existence of strong stationary solutions}

Strong stationary solutions are solutions of the stationary problem related to the one described in Section \ref{sec:1}.
More specifically, given $\alpha >0$, we say that a pair of functions $ (\u_\alpha, m_\alpha)$ is a strong stationary solution if the following conditions are met.

Functions $\u_\alpha$ and $m_\alpha$ must have at least the following regularity:
\[
 \u_\alpha \in \J_{0}\cap \H^2(\Omega)
 \quad \mbox{and} \quad
  m_\alpha \in  H^2(\Omega) ,
 \]
 
\noindent
and  $ (\u_\alpha, m_\alpha)$ must be a  solution of the weak formulation corresponding to the stationary system related to \eqref{sistema2}:
\begin{equation}
\label{StationarySystem}
\left\{\begin{array}{l} 
\displaystyle
- 2 \ {\rm div}\, (\nu(m_\alpha)D(\u_\alpha)) + \u_\alpha \cdot \nabla \u_\alpha + \nabla q_\alpha = -m_\alpha \chi +\f, 
\\
{\rm div}\, \u_\alpha = 0, 
\\
\displaystyle
 - \theta \Delta m_\alpha +\u_\alpha \cdot \nabla m_\alpha + U \displaystyle\frac{\partial m_\alpha}{\partial x_3} = 0, \qquad {\rm in} \; \Omega ,
\end{array}\right.
\end{equation}

\noindent
with boundary conditions related to \eqref{bc0}:
\begin{equation}
\label{bc0Stationary}
\left\{\begin{array}{l}
\u_\alpha = 0 \qquad {\rm on} \; S, 
\\
\u_\alpha \cdot \n = 0 \qquad {\rm on} \;  \Gamma,
\\
\nu(m_\alpha)[D(\u_\alpha)\n - \n\cdot(D(\u_\alpha)\n)\n] =\b_1  \qquad {\rm on} \;  \Gamma, 
\vspace{0.1cm}
\\
\theta \displaystyle\frac{\partial m_\alpha}{\partial \n} - Um_\alpha n_3 = 0
\qquad {\rm on} \;  \partial\Omega .
\end{array}\right.
\end{equation}

\noindent
Moreover, it is also required the following total mass condition related to \eqref{TotalMassConservation}:
\begin{equation}
\label{TotalMassConditionStationary}
 \int_\Omega  m_\alpha (x) dx = \alpha.
\end{equation}

In the above system, $q_\alpha$ is the hydrostatic pressure associated to the pair   $ (\u_\alpha, m_\alpha)$.

\vspace{0.2cm}
Existence of strong stationary solutions can be proved under certain conditions. 
For instance, in  Boldrini {\it et al.} \cite{BRMRM1} it is presented the following result,
which we restate here as a proposition for easy reference.
\begin{proposition} 
\label{ssolestac}
Assume that $\f \in \X$, $\b_1=0$ and also that $\nu$ a  strictly positive function of class  $C^1_b$.
Then, for $U$ sufficiently small, there exists a strong solution of \eqref{StationarySystem}, \eqref{bc0Stationary} and \eqref{TotalMassConditionStationary}.
\end{proposition}

\begin{remark}
When $b_1\neq 0$, it is possible to prove the existence of  weak stationary solutions.
\end{remark}

\subsection{Exponential $L^2$-stability of small strong stationary solutions}

Concerning strong stationary solutions, we prove the following stability result.
\begin{theorem}
\label{L2StabilityStationarySolutionResult}
Assume that $\f \in \X$ and that $\nu$ a   function of class  $C^1_b$ satisfying \eqref{nus}.
Then, any  small enough strong stationary solution   $ (\u_\alpha, m_\alpha)$ of \eqref{StationarySystem}, \eqref{bc0Stationary} and \eqref{TotalMassConditionStationary} is exponentially $L^2$-stable.

\end{theorem}

\begin{remark}
We stress that in this section $m_\alpha$ is not the same as the one in the previous sections.
\end{remark}

\subsection{Proof of Theorem \ref{L2StabilityStationarySolutionResult}}

Consider a particular  strong stationary solution $(\u_\alpha, m_\alpha)$ of 
\eqref{StationarySystem}, \eqref{bc0Stationary} and \eqref{TotalMassConditionStationary},
which according to Proposition \ref{ssolestac} exists at least for small enough $U$.

Consider also   an evolutive weak solution $(\u, m)$ and the following change of variables:
\[
\v=\u-\u_\alpha \quad{and} \quad \eta=m-m_\alpha .
\]

Then,  $(\v, \eta)$ must satisfy the weak formulation associated to the following system:
\begin{equation}
\label{sistema3}
  \left\{
  \begin{array}{l}
  \displaystyle\frac{\partial\v}{\partial t}-2\,{\rm div}
  (\nu(\eta+m_{\alpha})D(\v))-2\,{\rm div}
  (\nu(\eta+m_{\alpha})D(\u_{\alpha}))\\
  +2\,{\rm div}(\nu(m_{\alpha})D(\u_{\alpha}))+
  \v\cdot\nabla \v+\v\cdot\nabla \u_{\alpha}+
  \u_{\alpha}\cdot\nabla \v+\nabla(q-q_{\alpha})=-\eta\cdot \chi,\\
  {\rm div}\, \v=0,\\
  \displaystyle\frac{\partial\eta}{\partial t}-\theta\Delta\eta
  +\v\cdot\nabla\eta+\v\cdot\nabla
  m_{\alpha} +\u_{\alpha}\cdot\nabla\eta+U
  \displaystyle\frac{\partial\eta}{\partial x_3}=0,
  \ {\rm in}\ (0,T)\times\Omega ;
  \end{array}
  \right.
  \end{equation}
  
  \noindent
  with boundary conditions
  \begin{equation}
  \label{bcPerturbation}
  \left\{
  \begin{array}{l}
  \v=0 \ {\rm on}\ (0,T)\times S,
  \\
  \v\cdot\textbf{n} =0 \ {\rm on}\ (0,T)\times \Gamma,
  \\
  \nu(\eta+m_{\alpha})[D(\v)\n-\n\cdot(D(\v)\n)\n]
  \\
  \displaystyle
  \hspace{2cm}
  =
 -  (\nu (\eta + m_\alpha ) - \nu (m_\alpha) ) [D(\u_\alpha )\n]_{\mathbf{t}}

 ,
  \ {\rm on}\ (0,T)\times \Gamma,
  \\
  \theta\displaystyle\frac{\partial\eta}{\partial\n}-U\eta n_3 =0\ {\rm on}\ (0,T) ,
  \times\partial\Omega ,
  \end{array}
  \right.
  \end{equation}

 \noindent
and initial conditions given by
\begin{equation}
\label{ic3}
\v =\v_0, \quad \eta=\eta_0 \quad\mbox{at time} \; t=0.
\end{equation}

\noindent
Here, we used the notation
\[
[D(\u_\alpha )\n]_{\mathbf{t}} = D(\u_\alpha)\n-\n\cdot(D(\u_\alpha)\n)\n .
\]

\vspace{0.1cm}
In weak formulation, the  \eqref{sistema3} and \eqref{bcPerturbation} become: 
\begin{equation}
\label{12}
  \left\{
\begin{array}{l}
\displaystyle
 \left< \partial_t \v, \w\right> + 2(\nu(\eta+m_\alpha)D(\v), D(\w)) + 2(\nu(\eta+m_\alpha) D(\u_\alpha), D(\w)) 
  \\\noalign{\medskip}
  - 2 (\nu(m_\alpha) D(\u_\alpha), D(\w)) + (\v\cdot\nabla\v,\w) +  (\v\cdot\nabla\u_\alpha, \w) + (\u_\alpha\cdot\nabla \v, \w) =  
 \\
 - (\eta \cdot \chi, \w)
+ 2 \int_{\Gamma}   (\nu (\eta + m_\alpha ) - \nu (m_\alpha) ) [D(\u_{\alpha})\n]_{\mathbf{t}} \cdot \w~dS
\vspace{0.1cm}
  \\
  \displaystyle
 \left< \partial_t \eta, \phi \right> + \theta (\nabla\eta, \nabla\phi) + (\v\cdot\nabla\eta,\phi) + (\v\cdot\nabla m_\alpha, \phi) +   (\u_\alpha\cdot\nabla \eta, \phi) - U\left(\eta,
\frac{\partial\phi}{\partial x_3} \right) = 0,
\end{array}
\right.
\end{equation}
in $D'(0,T)$, for all $\w \in \J_0$ and $\phi \in B$.

\vspace{0.2cm}
To prove the result stated in Theorem \ref{L2StabilityStationarySolutionResult}, it is enough to prove the exponential decay of the perturbations $(\v, \eta)$.
For this, we need to obtain certain estimates for $\v$ and $\eta$.

We start by taking $\w = \v$ and $\phi = \eta$ in \eqref{12}; 
since $\displaystyle \left< \partial_t \v, \v \right> =  \frac12 \frac{d}{dt} |\v|^2$
and $  \displaystyle  \left< \partial_t \eta, \eta \right> = \frac12 \frac{d}{dt} |\eta|^2$,
we obtain that
\begin{equation}
\label{vwn}
\begin{array}{l}
\displaystyle
 \frac12 \frac{d}{dt} |\v|^2 + 2 (\nu(\eta+m_\alpha) D(\v), D(\v)) 
+ 2 \left( (\nu(\eta+m_\alpha) - \nu(m_\alpha)) D(\u_\alpha), D(\v) \right) 
\vspace{0.2cm}
\\
\displaystyle
 + (\v^n\cdot\nabla \u_\alpha, \v) + (\eta \cdot \chi, \v) 
 = 
 2 \int_{\Gamma}  (\nu (\eta + m_\alpha ) - \nu (m_\alpha) ) [D(\u_{\alpha})\n]_{\mathbf{t}} \cdot \v~dS ,
\end{array}
\end{equation}

\begin{equation}
\label{nuwn}
\displaystyle
\frac12 \frac{d}{dt} |\eta|^2 + \theta |\nabla\eta|^2 
+ (\v^n\cdot\nabla m_\alpha, \eta) - U\left(\eta, \frac{\partial \eta}{\partial x_3}\right) = 0.
\end{equation}

By using \eqref{nus}, H\"older and Young inequalities and also Sobolev embeddings, we obtain the following estimates.
From (\ref{vwn}), we have that
\begin{eqnarray*}
\frac12 \frac{d}{dt} |\v|^2 + 2 \nu_0 |D(\v)|^2 & \leq & 
2 \nu^\prime_1 |\eta|_4 |D(\u_\alpha)|_4 \ |D(\v)| 
+ |D(\u_\alpha)|_4 \ |\v|_4 \ |\v| + |\eta| \ |\v|, 
\\
&&  + 2   \nu_1 \|\eta \|_{L^4 (\Gamma)}  \| \nabla \u_{\alpha}\|_{L^4 (\Gamma)} \| \v \|_{L^2 (\Gamma)}  ,
\\
 & \leq & 
 2  C_1 (\nu^\prime_1 + 1)\|\nabla \u_\alpha\|_1 | \nabla \eta| \ |D(\v)|  
 \\
&& 
 + 2 \nu_1  \|\eta \|_{H^{1/2} (\Gamma)}     \| \nabla \u_{\alpha})\|_{H^{1/2} (\Gamma)} \| \v \|_{L^2 (\Gamma)}  
\\
& \leq & 2  C_1 (\nu^\prime_1 + 1)\|\nabla \u_\alpha\|_1 | \nabla \eta| \ |D(\v)|  
\\
&&
 + 2   \nu_1 | \nabla \eta | \,  \| \nabla  \u_{\alpha} \|_1 | D (\v ) |  
\\
& \leq &\frac{3}{2}\nu_0 |D(\v)|^2 +  C_{\nu_0} (4  C_1^2 (\nu^\prime_1 + 1)^2  ) + 4 \nu_1^2) \| \nabla \u_\alpha\|_1^2 \, | \nabla \eta|^2 .
\end{eqnarray*}

\noindent
This implies that
\begin{equation}
\label{ExponentialEstimatePreparatory1}
\frac{d}{dt} |\v|^2 +  \nu_0 |D(\v)|^2  \leq 
 C_{\nu_0} (4  C_1^2 (\nu^\prime_1 + 1)^2  ) + 4 \nu_1^2) \| \nabla \u_\alpha\|_1^2 | \nabla \eta|^2 .
\end{equation}

\vspace{0.1cm}
From (\ref{nuwn}), we get that
\begin{eqnarray*}
\frac12 \frac{d}{dt} |\eta|^2 + \theta |\nabla\eta|^2 & \leq &
|\nabla m_\alpha|_3 \ |\v|  \ |\eta|_6 + U|\eta| \ |\nabla\eta|
\\
 & \leq &
C_2 \| \nabla m_\alpha \|_1 \ | D(\v)|   |\nabla \eta| + U|\eta| \ |\nabla\eta|
\\
& \leq &\frac{\theta}{2}  |\nabla\eta|^2  + C_\theta \left(C_2^2 \| \nabla m_\alpha \|_1^2 + UC_0 \right) |\nabla\eta|^2 ,
\end{eqnarray*}

\noindent
which leads to
\begin{equation}
\label{ExponentialEstimatePreparatory2}
 \frac{d}{dt} |\eta|^2 + \theta |\nabla\eta|^2 
 \leq  
2 C_\theta \left(C_2^2 \| \nabla m_\alpha \|_1^2 + UC_0 \right) |\nabla\eta|^2 .
\end{equation}

By adding inequalities (\ref{ExponentialEstimatePreparatory1}) and (\ref{ExponentialEstimatePreparatory2}) , assuming tht
$\|\nabla \u_\alpha\|_1$, $\| \nabla m_\alpha \|_1$ and  $U$ small enough such that
\[
 C_{\nu_0} (4  C_1^2 (\nu^\prime_1 + 1)^2  ) + 4 \nu_1^2) \| \nabla \u_\alpha\|_1^2 
+ 2 C_\theta \left(C_2^2 \| \nabla m_\alpha \|_1^2 + UC_0 \right)
\leq \frac{\theta}{2} ,
\]

\noindent
we obtain that
\[
\frac{d}{dt}\left( |\v|^2 +  |\eta|^2 \right) +  \nu_0 |D(\v)|^2 + \frac{\theta}{2} |\nabla\eta|^2   \leq  0 .
\]

Since there is a constant $\overline{C}_{\nu_0\theta}$ such that 
$\overline{C}_{\nu_0\theta} \left( |\v|^2 +  |\eta|^2 \right) \leq  \nu_0 |D(\v)|^2 + \frac{\theta}{2} |\nabla\eta|^2 , $
the last inequality implies that 
\[
\frac{d}{dt}\left( |\v|^2 +  |\eta|^2 \right) + \overline{C}_{\nu_0\theta} \left( |\v|^2 +  |\eta|^2 \right)  \leq  0 ,
\]

\noindent
which gives the exponential decay of $ |\v(t)|^2 +  |\eta(t)|^2 $ as $t \rightarrow +\infty$, 
and so to the stated $L^2$-stability.
\hfill $\Box$


\end{document}